\documentclass[11pt,reqno]{amsart}
\usepackage{fullpage,times,graphicx,amssymb,amsmath,amsfonts,bbm,psfrag}
\usepackage{xcolor}
\usepackage[colorlinks,citecolor=bbluegray,linkcolor=ddarkbrown,urlcolor=blue,breaklinks]{hyperref}
\newtheorem{theorem}{Theorem}[section]
\newtheorem{proposition}[theorem]{Proposition}
\newtheorem{definition}[theorem]{Definition}

\newtheorem{corollary}[theorem]{Corollary}
\renewenvironment{proof}{\textbf{Proof.}}{\QED\bigskip}

\usepackage{algorithmic,algorithm}

\definecolor{ddarkbrown}{rgb}{0.5,0.2,0.05} \definecolor{bbluegray}{rgb}{0.05,0,0.5}

\newcommand{\BEAS}{\begin{eqnarray*}}
\newcommand{\EEAS}{\end{eqnarray*}}
\newcommand{\BEA}{\begin{eqnarray}}
\newcommand{\EEA}{\end{eqnarray}}
\newcommand{\BEQ}{\begin{equation}}
\newcommand{\EEQ}{\end{equation}}
\newcommand{\BIT}{\begin{itemize}}
\newcommand{\EIT}{\end{itemize}}
\newcommand{\BNUM}{\begin{enumerate}}
\newcommand{\ENUM}{\end{enumerate}}

\newcommand{\BA}{\begin{array}}
\newcommand{\EA}{\end{array}}

\newcommand{\ones}{\mathbf 1}

\newcommand{\reals}{{\mathbb R}}

\newcommand{\symm}{{\mbox{\bf S}}}  

\newcommand{\Span}{\mbox{\textrm{span}}}

\newcommand{\idm}{\mathbf{I}}

\newcommand{\QED}{~~\rule[-1pt]{6pt}{6pt}}

\newcommand{\argmin}{\mathop{\rm argmin}}

\usepackage[square]{natbib}

\usepackage{auto-pst-pdf} 

\let \oldsection \section
\renewcommand{\section}{\vspace{3ex plus 1ex}\oldsection}

\title{Regularized Nonlinear Acceleration\footnote{\textrm{A subset of these results appeared at the 2016 NIPS conference under the same title.}}}

\author{Damien Scieur}
\address{INRIA \& D.I.,\vskip 0ex
\'Ecole Normale Sup\'erieure, Paris, France.}
\email{damien.scieur@inria.fr}

\author{Alexandre d'Aspremont}
\address{CNRS \& D.I., UMR 8548,\vskip 0ex
\'Ecole Normale Sup\'erieure, Paris, France.}
\email{aspremon@ens.fr}

\author{Francis Bach}
\address{INRIA \& D.I.\vskip 0ex
\'Ecole Normale Sup\'erieure, Paris, France.}
\email{francis.bach@inria.fr}

\keywords{Acceleration, $\varepsilon$-algorithm, extrapolation.}
\date{\today}
\subjclass[2010]{}

\begin{document}

\begin{abstract}
We describe a convergence acceleration technique for unconstrained optimization problems. Our scheme computes estimates of the optimum from a nonlinear average of the iterates produced by any optimization method. The weights in this average are computed via a simple linear system, whose solution can be updated online. This acceleration scheme runs in parallel to the base algorithm, providing improved estimates of the solution on the fly, while the original optimization method is running. Numerical experiments are detailed on classical classification problems.
\end{abstract}

	\maketitle

\section{Introduction}\label{ss:intro}
Suppose we seek to solve the following optimization problem
\BEQ\label{eq:fprob}
\min_{x\in\reals^n} f(x)
\EEQ
in the variable $x\in\reals^n$, where $f(x)$ is strongly convex with parameter~$\mu$ with respect to the Euclidean norm, and has a Lipschitz continuous gradient with parameter $L$ with respect to the same norm. Assume we solve this problem using the \emph{fixed-point iteration}
\BEQ\label{eq:iter}
\tilde x_{i+1} = g( \tilde x_i), \quad \mbox{for $i=0,...,k$,}  \tag{FPI}
\EEQ
where $\tilde x_i\in \reals^n$ and $k$ is the number of iterations. This iteration is typically produced by an optimization algorithm, e.g. the gradient method with fixed step size, written
\BEQ\label{eq:grad-iter}
\tilde x_{i+1} = \tilde x_i -  h \nabla f (\tilde x_i), \quad \mbox{for $i=0,...,k$,}
\EEQ
with step length $h>0$. Here, we will focus on improving our estimates of the solution to problem~\eqref{eq:fprob} by tracking only the iterate sequence $\tilde x_i$ produced by an optimization algorithm, without any further calls to oracles on $g(x)$.

Since the publication of Nesterov's optimal first-order smooth convex minimization algorithm \citep{Nest83}, significant efforts have been focused on either providing more interpretable views on current acceleration techniques, or on replicating these complexity gains using different, more intuitive schemes. Early efforts sought to directly extend the original acceleration result in \citep{Nest83} to broader function classes \citep{Nemi85}, allow for generic metrics, line searches, produce simpler proofs \citep{Beck09,Nest03a} or adaptive accelerated algorithms \citep{Nest15}, etc. More recently however, several authors \citep{Dror14,Less16} have started using classical results from control theory to obtain numerical bounds on convergence rates that match the optimal rates. Others have studied the second order ODEs obtained as the limit for small step sizes of classical accelerated schemes, to better understand their convergence \citep{Su14,Wibi15}. Finally, recent results have also shown how to wrap classical algorithms in an outer optimization loop, to accelerate convergence and reach optimal complexity bounds \citep{Lin15} on certain structured problems.

Here, we take a significantly different approach to convergence acceleration stemming from classical results in numerical analysis. We use the iterates produced by any (converging) optimization algorithm, and estimate the solution directly from this sequence, assuming only some regularity conditions on the function to minimize. Our scheme is based on the idea behind Aitken's~$\Delta^2$-algorithm \citep{Aitk27}, generalized as the Shanks transform \citep{Shan55}, whose recursive formulation is known as the $\varepsilon$-algorithm \citep{Wynn56} (see e.g. \citep{Brez77,Sidi86} for a survey). In a nutshell, these methods fit geometrical models to linearly converging sequences, then extrapolate their limit from the fitted model. 

In a sense, this approach is more statistical in nature. It assumes an approximately linear model holds for iterations near the optimum, and estimates this model using the iterates. In fact, Wynn's algorithm \citep{Wynn56} is directly connected to the Levinson-Durbin algorithm \citep{Levi49,Durb60} used to solve Toeplitz systems recursively and fit autoregressive models (the Shanks transform solves Hankel systems, but this is essentially the same problem \citep{Hein11}). The key difference in these extrapolation techniques is that estimating the autocovariance operator $A$ is not required, as we only focus on the limit. Moreover, the method presents strong links with the conjugate gradient when applied to unconstrained quadratic optimization, but does not further calls to the operator.

We start from a formulation of these techniques known as Anderson Acceleration \citep{anderson1965iterative}, Me{\v{s}}ina's Algorithm \citep{Mesi77} or minimal polynomial extrapolation (MPE) \citep{Sidi86,Smit87}. They use the minimal polynomial of the linear operator driving iterations to estimate the optimum by a nonlinear average of the iterates (i.e. computing a weighted average using weights which are nonlinear functions of the iterates). 

Our contribution here is to regularize this procedure and produce explicit bounds on the distance to optimality by controlling stability, thus explicitly quantifying acceleration. We show that these extrapolation algorithms reach optimal performance (asymptotically) and describe several numerical examples where these stabilized estimates often speed up convergence by an order of magnitude. So far, for all the techniques cited above, no proofs of convergence of the estimates were given when the estimation process became unstable. Furthermore, the acceleration scheme runs in parallel with the original algorithm, providing improved estimates of the solution on the fly, while the original method is progressing, so its numerical complexity is marginal.

The paper is organized as follows. In Section \ref{sec:acc_fixedpoint} we recall basic results behind the acceleration for linear iterations. Then, in Section \ref{sec:reg_acc}, we generalize these results to nonlinear iterations and show  how to fully control the impact of nonlinearity. We use these results to derive explicit bounds on the acceleration performance of our estimates. In Section \ref{sec:ext} we connect the acceleration methods to the conjugate gradient method and Nesterov's method. Finally, we present numerical results in Section~\ref{sec:numres}.

\section{Convergence Acceleration}\label{sec:acc_fixedpoint}
We begin by recalling the core arguments behind convergence acceleration. These ideas have taken various forms over time, known for example as {\em Anderson acceleration} \citep{anderson1965iterative}, the \textit{Eddy-Mesina method} \citep{Mesi77,Eddy79} and {\em minimal polynomial extrapolation} \citep{Caba76,Smit87}. The core idea behind these methods is to use a Taylor expansion of the function $g$ in~\eqref{eq:iter} to approximate the fixed point iterations by a vector autoregressive model, then compute a weighted mean of the iterates $\tilde x_i$ to produce a better estimate of the limit $x^*$. In this paper, we assume $x^*$ unique.

Suppose $g(x)$ is differentiable and let $G$ be the Jacobian of $g$ evaluated at $x^*$. In the rest of the paper, we assume $G$ to be symmetric, positive semi-definite and $G \preceq \sigma I$, with $\sigma <1$. Equation \eqref{eq:iter} becomes
\[
	\tilde x_{i+1} = g(x^*) + G(\tilde x_i-x^*) + O(\|\tilde x_i-x^*\|^2), \quad \mbox{for $i=1,\ldots,k.$}
\]
By neglecting the second order term, and because $g(x^*) = x^*$, we obtain the \emph{linear fixed-point iteration}
\BEQ
	x_{i+1}-x^* = G(x_{i}-x^*),   \label{eq:linear_fpi} \tag{LFPI}
\EEQ
where $x_0 = \tilde x_0$ and we recognize here a vector autoregressive process. When using the fixed-step gradient method in~\eqref{eq:grad-iter} for example, if $\nabla^2f$ is the Hessian matrix of $f(x)$, we get
\[
	x_{i+1}-x^* = \underbrace{(\idm-\nabla^2f(x^*))}_{= G}(x_{i}-x^*).
\]
Because $\|G\|_2 \leq \sigma <1$, the iterates $x_k$ converges to $x^*$ at a linear rate, with
\[
	\|x_i-x^*\| \leq \sigma\|x_{i-1}-x^*\| \leq \sigma^i\|x_0-x^*\|,
\]
where $\|\cdot\|$ stands for the Euclidean norm here and throughout the paper. We will now see how to improve convergence rates using a linear combination of the previous iterates.

Suppose we run $k$ iterations of \eqref{eq:linear_fpi}, a linear combination of iterates $x_i$ with coefficients $c_i$ reads
\BEA
	\sum_{i=0}^k c_i x_{i} & = & \sum_{i=0}^k c_i x^* + \sum_{i=0}^k c_i G(x_i-x^*) \nonumber\\
	& = & \left(\sum_{i=0}^k c_i\right) x^* + \left(\sum_{i=0}^k c_i G^i\right)(x_0-x^*). \label{eq:linear_combination}
\EEA
Now define the polynomial 
\BEA
	p(z) \triangleq \sum_{i=0}^k c_i z^i, \label{eq:polynomial}
\EEA
we can write \eqref{eq:linear_combination} more concisely in terms of the matrix polynomial $p(G)$, setting $p(1) = \sum_{i=0}^k c_i = 1$ without loss of generality, to get
\[
	\sum_{i=0}^k c_i x_{i} = x^* + \underbrace{p(G) (x_0-x^*).}_{\text{Error term}}
\]
Ideally, we need to find $c$ (or equivalently $p$) which minimizes the error term $p(G) (x_0-x^*)$. We will study the error when linearly combining the last $k+1$ iterates $x_i$, assuming we have an algorithm computing this optimal combination, i.e.
\[
\left\| \sum_{i=0}^k c^\star_i x_{i} - x^* \right\| = \min_{\{c \in\reals^{k+1}:\,c^T\ones = 1\}} \left\|\sum_{i=0}^k c_iG^i(x_0-x^*)\right\| =\min_{\{p \in \reals_k[x]:\,p(1) = 1\}} \left\|p(G)(x_0-x^*)\right\|
\]
where $\reals_k[x]$ is the subspace of polynomials of degree at most $k$ and 
\[
c^\star = \argmin_{\{c \in\reals^{k+1}:\,c^T\ones = 1\}} \left\|\sum_{i=0}^k c_iG^i(x_0-x^*)\right\|.
\]
The next proposition produces an uniform bound on the value of this error using Chebyshev polynomials.

\begin{proposition} \label{prop:rateconv_linearacc}
Suppose the iterates $x_i$ for $i=0,\ldots,k$ are computed using \eqref{eq:linear_fpi}, with $G$ the Jacobian of $g$, assumed to be symmetric, satisfying $0 \preceq G \preceq \sigma I$ for $\sigma <1$. Let $x^*$ be the fixed point of $g$. The $\ell_2$ norm of the error is bounded, with
	\BEA
		\left\| \sum_{i=0}^k c^\star_i x_{i} - x^* \right\| \leq 
		\begin{cases}
		\dfrac{2\beta^k}{1+\beta^{2k}} \, \|x_0-x^*\| \qquad & \text{if $k < m$} \\
		0 & \text{otherwise}\\
		\end{cases} \label{eq:rate_linear}
	\EEA
	where $m$ is the number of distinct eigenvalues of $G$ and
	\BEQ \label{eq:defbeta}
		\beta = \frac{1-\sqrt{1-\sigma}}{1+\sqrt{1-\sigma}} < 1.
	\EEQ
\end{proposition}
\begin{proof}
Because $G$ is symmetric, it admits the eigenvalue decomposition
	\[
		G = Q^*\Lambda Q,
	\]
	where $\Lambda$ is the diagonal matrix of eigenvalues $\{\lambda_i, \; i=1, \ldots , m\}$, and $Q$ is a unitary matrix.  For any $p\in\reals_k[x]$, we have
	\BEAS
\|p(G)(x_0-x^*)\|_2   & = & \|Q^*p(\Lambda)Q(x_0-x^*)\|_2  \\
& \leq & \|p(\Lambda) \|_2~\|(x_0-x^*)\|_2 \nonumber \\
& = & \max_{i=1,\ldots,m} |p(\lambda_i)|~\|(x_0-x^*)\|_2 . 
	\EEAS
	
	First, assume $k \geq m$. The Cayley-Hamilton theorem means that if $p(x)$ is the characteristic polynomial of $G$, then $p(G)=p(\Lambda)=0$. By assumption none of the $\lambda_i$ is equal to 1 (we assumed $\lambda_i \in [0,\sigma]$ with $\sigma < 1$) so we can normalize $p$ so that $p(1)=1$ and $p(\lambda_i) = 0$ for all $i =1,\ldots,m$.

We now assume $k < m$. We have, for any polynomial $p$,
\[
\max_{i=1,\ldots,m} |p(\lambda_i)| \leq  \max_{\lambda \in [\lambda_{\min},\lambda_{\max}]} |p(\lambda)|
\]
Because $0 \preceq G \preceq \sigma$, we have $0 \leq \lambda_i \leq \sigma$, so the error bound becomes
\BEQ \label{eq:problem_minimax_cheby}
\min_{\{p \in \reals_k[x]:\, p(1) = 1\}} \|p(G)(x_0-x^*)\|_2 \leq  \min_{\{p \in \reals_k[x]:\, p(1) = 1\}}~\max_{\lambda \in [0,\sigma]} |p(\lambda)|~\| x_0-x^* \|_2
\EEQ
where the right hand side involves a minmax problem, explicitly solved using Chebyshev's polynomials. Let~$C_k$ be the Chebyshev polynomial of degree $k$. By definition, $C_k$ is a monic polynomial (i.e. a polynomial whose leading coefficient is one) solving
	\[ 
	C_k (x) \triangleq \argmin_{\{p\in\reals_k[x]:\, p(1)=1\}} \;\max_{x\in[-1,1]}|p(x)|.
	\]
\citet{Golu61} use a variant of $C_k(x)$ to solve the problem in \eqref{eq:problem_minimax_cheby}, whose solution is a rescaled Chebyshev polynomial given by 
	\BEQ\label{eq:rescaledcheby}
	T_{k}(x,\sigma) = \frac{C_k(t(x,\sigma))}{C_k(t(1,\sigma))}, \quad \mbox{where}\quad t(x,\sigma) =\frac{2x-\sigma}{\sigma},
	\EEQ
where $t(x,\sigma)$ is simply a linear mapping from interval $[0,\sigma]$ to $[-1,1]$. Moreover, they show
	\BEQ\label{eq:chebysol}
		\min_{\{p\in\reals_k[x]:\, p(1) = 1\}}~\max_{\lambda \in [0,\sigma]}|p(\lambda)| = \max_{\lambda \in [0,\sigma]}|T_{k}(\lambda,\sigma)| = |T_{k}(\sigma,\sigma)| = \frac{2\beta^k}{1+\beta^{2k}},
	\EEQ
	where $\beta$ is given by
	\[
	\beta = \frac{1-\sqrt{1-\sigma}}{1+\sqrt{1-\sigma}} < \sigma <1.
	\]
	Injecting the result of \eqref{eq:chebysol} in \eqref{eq:problem_minimax_cheby} yields the desired result.
\end{proof}

\begin{figure}[ht]
\begin{center}
			\psfrag{x}[c][b]{$x$}
			\psfrag{Txs}[b][t]{$T_k(x,\sigma)$}
			\psfrag{sigm}[l][r]{$\sigma$}
			\includegraphics[width=0.5\textwidth]{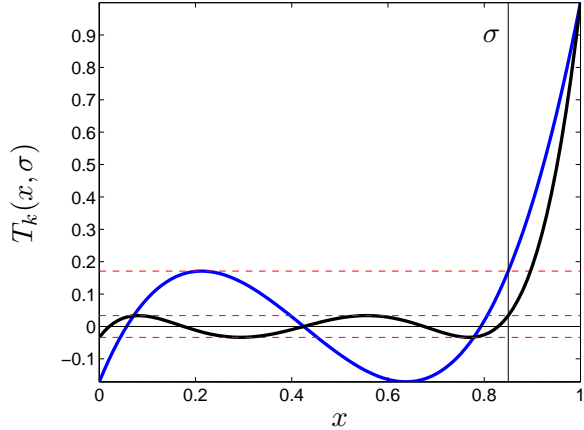}
\caption{We plot both $T_{3}(x,\sigma)$ (blue) and $T_{5}(x,\sigma)$ (black) for $x\in[0,1]$ and $\sigma=0.85$. The maximum value of the image of $[0,\sigma]$ by $T_k$ is clearly smaller when $k$ grows, implying a better rate of convergence. \label{fig:cheb}}
\end{center}
\end{figure}

\begin{corollary}
In the case of the gradient method applied on quadratic function with eigenvalues bounded in the interval $[\mu,L]$ (this correspond also to a $L$-smooth and $\mu$-strongly convex function), we have
\[
	\sigma = 1-\frac{\mu}{L} <1.
\]
By consequence, the bound \eqref{eq:rate_linear} becomes
	\BEAS
		\left\| \sum_{i=0}^k c^\star_i x_{i} - x^* \right\| \leq 
		\begin{cases}
		\dfrac{2\beta^k}{1+\beta^{2k}} \, \|x_0-x^*\| \qquad & \text{if $k < m$} \\
		0 & \text{otherwise}\\
		\end{cases} 
	\EEAS
	with
	\[
		\beta = \frac{1-\sqrt{\mu/L}}{1+\sqrt{\mu/L}}.
	\]
\end{corollary}

The proof of this proposition suggests $T_{k}(x,\sigma)$ is a good universal solution for the convergence acceleration problem and we plot both $T_{3}(x,\sigma)$ and $T_{5}(x,\sigma)$ for $x\in[0,1]$ and $\sigma=0.85$ in Figure~\ref{fig:cheb}. This solution is called the {\em Chebyshev semi-iterative method} in \citep{Golu61} and was further studied by e.g. \citep{Nemi84}. Combining the $k$ iterates $x_i$ using the coefficients $c_i$ in $T_{k}(x,\sigma)$, we ensure 
\[
	\left\|\sum_{i=0}^k c_i x_i - x^* \right\| \lesssim (1- \sqrt{1-\sigma})^k \|x_0-x^*\| \ll \sigma^k \|x_0-x^*\|
\]
which means convergence is indeed accelerated. However, this method has some key drawbacks. First, we need to know $\sigma$ to form $T_{k}(x,\sigma)$, which is not always the case. For example, in the case of the gradient method, $\sigma$ depends on the smoothness constant $L$ and the strong convexity constant $\mu$. In the general non-linear case, $\sigma$ depends on the spectrum of the Jacobian at the optimum, which is clearly not observed. Second, the algorithm does not allow us to control the magnitude of the coefficients in the polynomial $T(x,\sigma)$, which has a strong impact on the stability of the algorithm in the presence of numerical errors, or when the iterates are generated by a non-linear function $g$.

Because of stability issues with Chebyshev acceleration, we focus now on a method which will approximately minimize the error $\|p(G)(x_0-x^*)\|_2$. Since we of course do not observe $G$ and $x^*$ we will work with the residuals
\BEQ
	\tilde r_i = \tilde x_{i+1} - \tilde x_i = g(\tilde x_i)-\tilde x_i \label{eq:residual_nonlinear},
\EEQ 
when $g$ is a linear function \eqref{eq:linear_fpi} this becomes
\BEA
	r_i = x_{i+1} - x_{i} = (G-I)(x_i-x^*). \label{eq:residual_linear}
\EEA
A linear combination of residuals $r_i$ with coefficients $c_i$ is written
\[
	\sum_{i=0}^k c_ir_i = (G-I)\sum_{i=0}^kc_i(x_i-x^*) = (G-I)p(G)(x_0-x^*).
\]
We recognize the error term we wanted to minimize, multiplied by the matrix $(G-I)$. Using the coefficients which minimize this alternative quantity will approximately minimize the error, as stated in the following proposition.
\begin{proposition} \label{prop:linear_rate_acceleration}
	Let $p^*(x)$ be the polynomial solving
	\BEAS
		p^*(x) = \argmin_{\{p\in\reals_k[x] :\, p(1) = 1\}} \|(G-I)p(G)(x_0-x^*)\|_2, \label{eq:min_polynomial}
	\EEAS
	whose coefficients, written $c^*$, satisfy
	\BEQ\label{eq:c-rre}
		c^* = \argmin_{\{c \in \reals^{k+1} :\, c^T\ones = 1\}} \left\|\sum_{i=0}^k c_i r_i \right\|_2.
	\EEQ
The iterates $x_i$ defined in \eqref{eq:linear_fpi} averaged with coefficients $c^*$ satisfy
	\BEQ \label{eq:error_extrapolation}
		\left\|\sum_{i=0}^k c^*_i x_i-x^*\right\| \leq \frac{1}{1-\sigma}~\min_{\{c \in\reals^{k+1}:\,c^T\ones = 1\}} \left\|\sum_{i=0}^k c_iG^i(x_0-x^*)\right\|,
	\EEQ
where we have assumed $0\preceq G \preceq \sigma I$, with $\sigma <1$.
\end{proposition}
\begin{proof}
	By definition of $c^*$, and using \eqref{eq:linear_fpi},
	\BEAS
		\left\|\sum_{i=0}^k c^*_i x_i-x^*\right\| & = & \|p^*(G)(x_0-x^*)\|, \\
		& = & \|(G-I)^{-1}(G-I)p^*(G)(x_0-x^*)\|, \\
		& \leq & \|(G-I)^{-1}\|~\|(G-I)p^*(G)(x_0-x^*)\|.
	\EEAS
	By using the definition of $p^*$,
	\[
		 \|(G-I)p^*(G)(x_0-x^*)\| = \min_{\{p\in\reals_k[x] :\, p(1) = 1\}} \|(G-I)p(G)(x_0-x^*)\|_2.
	\]
	We can bound this last error term because $\|G-I\| \leq 1$,
	\BEAS
		\min_{\{p\in\reals_k[x] :\, p(1) = 1\}} \|(G-I)p(G)(x_0-x^*)\| & \leq & \min_{\{p\in\reals_k[x] :\, p(1) = 1\}} \|(G-I)\| \|p(G)(x_0-x^*)\|, \\
		& \leq &  \min_{\{p\in\reals_k[x] : p(1) = 1\}} \|p(G)(x_0-x^*)\|.
	\EEAS
Using the fact that $ \|(G-I)^{-1}\|\leq \frac{1}{1-\sigma}$ yields the desired result.
\end{proof}

This leads to the following acceleration algorithm.
\begin{algorithm}[H]
\caption{Nonlinear Acceleration of Convergence}
\label{algo:acc_fixedpoint}
\begin{algorithmic}[1]
\REQUIRE Iterates $x_0,x_1,\ldots,x_{k+1} \in \reals^d$.
\STATE Form $R = [r_0,...,r_k]$
\STATE Solve 
\[
c^* = \argmin\limits_{\{c \in\reals^{k+1}:\, c^T1 = 1\}} \|Rc\|
\]
\ENSURE Approximation of  $x^*$ ensuring \eqref{eq:error_extrapolation}, computed as $\sum_{i=0}^{k} c_i^* x_i$
\end{algorithmic}
\end{algorithm}

This acceleration algorithm is called nonlinear because the coefficients $c_i$ vary with of $x_i$. This method is also known as Anderson acceleration \citep{anderson1965iterative}, the Eddy-Mesina algorithm \citep{Mesi77,Eddy79}, Minimal Polynomial Extrapolation \citep{Caba76}, or Reduced Rank Extrapolation \citep{Sidi86,Smit87}. There are small variations between all these methods, which lie in the way they solve the minimization problem in~\eqref{eq:c-rre}. The next proposition gives us an explicit solution, involving the inversion of a $k$-by-$k$ matrix.

\begin{proposition} \label{prop:explicit_sol_linear}
The explicit solution of the problem
\BEA
	 c^* = \argmin\limits_{c^T1 = 1} \|Rc\| \label{eq:argmin_cstar}
\EEA
in the variable $c\in\reals^k$, where $R$ is a $d\times k$ matrix assumed to be of rank $k$, is given by
\BEA
	c^* = \frac{(R^TR)^{-1}\ones}{\ones^T(R^TR)^{-1}\ones}. \label{eq:linearsystem_cstar}
\EEA
\end{proposition}
\begin{proof}
Let $\mu$ be the dual variable of the equality constraint. Both $c^*$ and $\mu^*$ should satisfy the KKT system
\BEA \label{eq:kkt_system}
\begin{bmatrix}
2R^TR & \ones \\ \ones^T & 0
\end{bmatrix}
\begin{pmatrix}
c^* \\ \mu^*
\end{pmatrix}
=
\begin{pmatrix}
0 \\ 1
\end{pmatrix}
\EEA
This block matrix can be inverted explicitly, with
\[
\begin{bmatrix}
2R^TR & \ones \\ \ones^T & 0
\end{bmatrix}^{-1}
=
\frac{1}{\ones^T (R^T R)^{-1} \ones}\begin{bmatrix}
\frac{1}{2}(R^T R)^{-1}\left(\ones^T (R^T R)^{-1} \ones I - \ones\ones^T (R^T R)^{-1}\right) & (R^TR)^{-1}\ones \\
\ones^T(R^TR)^{-1} & -2
\end{bmatrix}.
\]
Using this inverse we easily solve the linear system, which gives the result in~\eqref{eq:linearsystem_cstar}.
\end{proof}

In practice of course, instead of computing the inverse of the matrix $R^TR$, we solve the linear system
\[
	R^TR z = \ones,
\]
then set $c^* = z/(\ones^Tz)$. This formula is used in Anderson Acceleration algorithm and Mesina method. Other algorithms usually force the coefficient $c_k$ to be equal to one, solve the remaining linear system, then normalize the vector. However, these alternative strategies are harder to analyze when the iterates are generated by a non-linear function $g$. We will now apply this acceleration algorithm on gradient method for nonlinear functions and compute its rate of convergence.

\section{Regularized Nonlinear Acceleration of Convergence}\label{sec:reg_acc}
So far, we have only considered linear functions $g$ in~\eqref{eq:linear_fpi}, without perturbations, when computing the iterates $x_i$. In general, the fixed-point iteration~\eqref{eq:iter} is usually generated by a nonlinear function $g$, thus inducing a second order error term in $O(\|x_i-x^*\|^2)$ compared to the dynamics in \eqref{eq:linear_fpi}.

Here, in \S\ref{ssec:sensitivity} we first give a bound on the deviation error when there are perturbations in \eqref{eq:linear_fpi}. In \S\ref{ssec:regularized_algo} we then derive a regularized version of Algorithm~\ref{algo:acc_fixedpoint} which better controls the impact of perturbations. We then study the impact of regularization on the solution when there are no perturbations in \S\ref{ssec:regularized_cheby}. Finally, in \S\ref{ssec_rate_conv_regalgo} we gather the results of the previous sections to bound the rate of convergence of the regularized acceleration algorithm.

\subsection{Sensitivity Analysis} \label{ssec:sensitivity}
We now study the sensitivity of the acceleration algorithm to perturbations. Consider the following perturbed linear fixed point iteration
\BEQ \label{eq:perturb_LFPI}
	\tilde x_{i+1} -x^* = g(\tilde x_i) - x^* \nonumber  = G(\tilde x_i - x^*) + e_i \tag{Pert. LFPI}
\EEQ
where $e_i$ is the noise injected in $x_{i+1}$ at iteration $i$. For now, we do not assume any structure on the noise, so $e_i$ may be the nonlinearity of $g$, stochastic noise, roundoff error, \emph{etc}. The iterates of this process will be compared to their noiseless counterpart,
\[
	x_{i+1}-x^* =  G(x_i - x^*),
\]
with $x_0 = \tilde x_0$. We now apply our acceleration algorithm on the sequences $x_i$ and $\tilde x_i$ and compare the results. We first form the residuals,
\[
	r_i =  g( x_{i})- x_i = x_{i+1} - x_i 
\qquad\mbox{and}\qquad
	\tilde r_i  =  g(\tilde x_{i})-\tilde x_i = \tilde x_{i+1} - \tilde x_i.
\]
Consider the matrices of residuals $R = [r_0,\ldots,r_k]$ and $\tilde R = [\tilde r_0,\ldots,\tilde r_k]$. We write $P$ the \emph{perturbation matrix} defined as
\BEA \label{eq:def_P}
	P \triangleq \tilde R^T \tilde R - R^TR.
\EEA
The next proposition describes the sensitivity of Algorithm \ref{algo:acc_fixedpoint} using $R$ and $P$.
\begin{proposition}\label{prop:sensitivity_acc}
	Let the sequences $x_i$ be generated by \eqref{eq:linear_fpi} and $\tilde x_i$ by \eqref{eq:perturb_LFPI}, with $x_0 = \tilde x_0$, with $R$ and $\tilde R$ the residual matrices defined above and $P$ the perturbation matrix in \eqref{eq:def_P}. Assume $c^*$ and $\tilde c^*$ are computed using formula \eqref{eq:linearsystem_cstar} with matrices $R$ and $\tilde R$ respectively. Let
	\BEA
		\Delta \tilde c^* \triangleq \tilde c^* - c^* \label{eq:def_delta_c}.
	\EEA
	Then the norm of $\Delta \tilde c^*$ is bounded by
	\BEQ
		\|\Delta \tilde c^*\| \leq \|P\| \|(R^T R + P)^{-1}\| \|c^*\|. \label{eq:bound_delta_c}
	\EEQ
\end{proposition}
\begin{proof}
We start with the sequence $\tilde x_i$. Let $\tilde\mu^*$ be the dual variable of the equality constraint of \eqref{eq:argmin_cstar}. Both $\tilde c^* = c^* + \Delta \tilde c^*$ and $\tilde \mu^* = \mu^* + \Delta \mu^*$ should satisfy the KKT system
\[
	\begin{bmatrix}
		2\tilde R^T \tilde R & \ones \\ \ones^T & 0
	\end{bmatrix}
	\begin{pmatrix}
		\tilde c^* \\ \tilde \mu^*
	\end{pmatrix}
	=
	\begin{pmatrix}
		0 \\ 1
	\end{pmatrix}
	\quad
	\Leftrightarrow
	\quad
	\begin{bmatrix}
		2(R^TR + P) & \ones \\ \ones^T & 0
	\end{bmatrix}
	\begin{pmatrix}
		c^* + \Delta \tilde c^* \\ \mu^* + \Delta \mu^*
	\end{pmatrix}
	=
	\begin{pmatrix}
		0 \\ 1
	\end{pmatrix}.
\]
Indeed, using the definition of $c^*$ and $\mu^*$ in \eqref{eq:kkt_system},
\BEAS
	\begin{bmatrix}
		2(R^TR + P) & \ones \\ \ones^T & 0
	\end{bmatrix}
	\begin{pmatrix}
		c^* + \Delta \tilde c^* \\ \mu^* + \Delta \mu^*
	\end{pmatrix} 
	& = & 
	\begin{bmatrix}
		2R^TR & \ones \\ \ones^T & 0
	\end{bmatrix}
	\begin{pmatrix}
		c^* \\ \mu^*
	\end{pmatrix}
	+
	\begin{bmatrix}
		2R^TR & \ones \\ \ones^T & 0
	\end{bmatrix}
	\begin{pmatrix}
		\Delta \tilde c^* \\ \Delta \mu^*
	\end{pmatrix}
	+
	\begin{bmatrix}
		2P & 0 \\ 0 & 0
	\end{bmatrix}
	\begin{pmatrix}
		c^* + \Delta \tilde c^* \\ \mu^* + \Delta \mu^*
	\end{pmatrix},
	\\
	& = & \begin{pmatrix}
	0 \\ 1
	\end{pmatrix}
	+
	\begin{bmatrix}
		2R^TR & \ones \\ \ones^T & 0
	\end{bmatrix}
	\begin{pmatrix}
		\Delta \tilde c^* \\ \Delta \mu^*
	\end{pmatrix}
	+
	\begin{bmatrix}
		2P & 0 \\ 0 & 0
	\end{bmatrix}
	\begin{pmatrix}
		c^* + \Delta \tilde c^* \\ \mu^* + \Delta \mu^*
	\end{pmatrix}.
 \EEAS
 With this simplification, the system becomes
 \[
 	\begin{bmatrix}
 		2R^TR & \ones \\ \ones^T & 0
 	\end{bmatrix}
 	\begin{pmatrix}
 		\Delta \tilde c^* \\ \Delta \mu^*
 	\end{pmatrix}
 	+
 	\begin{bmatrix}
 		2P & 0 \\ 0 & 0
 	\end{bmatrix}
 	\begin{pmatrix}
 		c^* + \Delta \tilde c^* \\ \mu^* + \Delta \mu^*
 	\end{pmatrix}
 	=
 	\begin{pmatrix}
 		0 \\ 0
 	\end{pmatrix}.
 \]
 It remains to isolate $c^*$,\[
  	\begin{bmatrix}
  		2(R^TR + P)  & \ones \\ \ones^T & 0
  	\end{bmatrix}
  	\begin{pmatrix}
  		\Delta \tilde c^* \\ \Delta \mu^*
  	\end{pmatrix}
  	=
  	\begin{pmatrix}
  		2Pc^* \\ 0
  	\end{pmatrix}.
  \]
  The explicit solution in obtained by inverting the block matrix, and is written
  \[
  	\Delta \tilde c^* = \left(I-\frac{(R^TR + P)^{-1}\ones \ones^T}{\ones^T (R^TR + P)^{-1} \ones}\right) (R^TR + P)^{-1}Pc^* .
  \]
  We can bound the norm of $\Delta \tilde c^*$ by
\[
   	\|\Delta \tilde c^*\| = \left\|I-\frac{(R^TR + P)^{-1}\ones \ones^T}{\ones^T (R^TR + P)^{-1} \ones}\right\| \|(R^TR + P)^{-1}\|\|P\|\|c^*\|.
\]
Because the first factor is the norm of a projector of rank $k-1$, its value is bounded by $1$, so we get the desired result.
\end{proof}

This proposition bounds the relative error on $\tilde c^*$ in comparison with $c^*$. We will see that the perturbation magnitude can be arbitrarily large, which is the key issue with the convergence results in \citep[\S7]{Smit87}. Even when $\|P\|$ is small, the term $\|(R^TR + P)^{-1}\|$ is problematic. Our problem here is the structure of the residuals matrix $R$,
\[
	R = [r_0,Gr_0,G^2r_0,...,G^kr_0],
\]
which matches exactly the structure of \emph{Krylov matrices}, i.e. square matrices $K$ formed using a matrix $M$ and a vector $v$, and computed as $K=[v,Mv,M^2v,..., M^{k}v]$. \citet{Tyrt94} showed that the condition number of Krylov matrices (see Section \ref{ssec:conjgrad}) is lower bounded by a function which grows exponentially with $k$. Now, the error bound \eqref{eq:bound_delta_c} contains the norm of the inverse of a perturbed squared Krylov matrix, which makes the situation even worse. In other words, even if the perturbations are small, their impact on the solution can be arbitrarily large. Even in practical cases where $k$ is small (for example, $k = 5$), $\tilde R^T \tilde R$ is usually a singular or nearly-singular matrix. This particular issue means the linear system $(R^TR)^{-1}\ones$ in \eqref{eq:linearsystem_cstar} needs to be regularized.

\subsection{Regularized Nonlinear Acceleration of Convergence} \label{ssec:regularized_algo}
In this section, we will analyze the following acceleration algorithm, which uses Tikhonov regularization to solve the linear system in~\eqref{eq:linearsystem_cstar}.
\begin{algorithm}[H]
\caption{Regularized Nonlinear Acceleration (RNA)}
\label{algo:acc_fixedpoint_reg}
\begin{algorithmic}[1]
\REQUIRE Iterates $\tilde x_0, \tilde x_1,..., \tilde x_{k+1}\in\reals^d$ produced by~\eqref{eq:iter}, and a regularization parameter $\lambda > 0$.
\STATE Compute $\tilde R = [\tilde r_0,...,\tilde r_k]$, where $\tilde r_i = \tilde x_{i+1}-\tilde x_i$
\STATE Solve 
\[
\tilde c^*_\lambda = \argmin\limits_{c^T1 = 1} \|\tilde Rc\|^2 + \lambda \|c\|^2,
\] 
or equivalently solve $(\tilde R^T\tilde R+\lambda I) z = \ones$ then set $\tilde c^*_\lambda = {z}/{\ones^Tz}$.
\ENSURE Approximation of $x^*$ computed as $\sum_{i=0}^{k} ( \tilde c_{\lambda}^*)_i \tilde x_i$
\end{algorithmic}
\end{algorithm}

Regularization controls the norm of the coefficients produced by the algorithm and reduces the impact of perturbations, as shown in the following proposition.
\begin{proposition}\label{prop:stability_reg_algo}
Consider the sequences $x_i$ satisfying \eqref{eq:linear_fpi} and $\tilde x_i$ satisfying \eqref{eq:perturb_LFPI} with $x_0 = \tilde x_0$. Let $c_\lambda^*$ and $\tilde c_\lambda^*$ the output of Algorithm \ref{algo:acc_fixedpoint_reg} with parameter $\lambda$ applied to $x_i$ and $\tilde x_i$ respectively. Let $R$ and $\tilde R$ the matrices of residuals and $P$ be defined in \eqref{eq:def_P}. Define $\Delta \tilde c^*_\lambda = \tilde c^*_\lambda - c_\lambda^*$. Then, we have the following bounds,
\BEA
	\|\tilde c^*_{\lambda}\| & \leq & \sqrt{\frac{\lambda + \|\tilde R\|^2}{(k+1)\lambda}},\\
	\|\Delta \tilde c_\lambda^*\| & \leq & \frac{\|P\|}{\lambda}\|c_\lambda^*\|,
\EEA
which control the stability of the solution $\tilde c_\lambda^*$.
\end{proposition}
\begin{proof}
Using the same proof technique of Propositions \ref{prop:explicit_sol_linear} and \ref{prop:sensitivity_acc}, we have
\BEA
	\tilde c^*_\lambda & = &  \frac{(\tilde R^T \tilde R +\lambda I)^{-1}\ones}{\ones^T(\tilde R^T \tilde R + \lambda I)^{-1}\ones}, \label{eq:sol_ctildelambda} \\
	\Delta \tilde c^*_\lambda & = & \left(I-\frac{(\tilde R^T \tilde R +\lambda I)^{-1}\ones \ones^T}{\ones^T (\tilde R^T \tilde R +\lambda I)^{-1} \ones}\right)   (\tilde R^T \tilde R +\lambda I)^{-1}Pc^*_\lambda. \label{eq:sol_deltactildelambda}
\EEA
We begin by the bound on $\tilde c^*_\lambda$. Indeed, with \eqref{eq:sol_ctildelambda},
\BEAS
	\|\tilde c^*_\lambda \|^2 & = & \frac{\ones^T(\tilde R^T \tilde R +\lambda I)^{-2}\ones}{(\ones^T(\tilde R^T \tilde R + \lambda I)^{-1}\ones)^2}, \\
	& \leq & \frac{1}{k+1}\max_{\|v\| = 1} \frac{v^T(\tilde R^T \tilde R +\lambda I)^{-2}v}{(v^T(\tilde R^T \tilde R + \lambda I)^{-1}v)^2},\\
	& = & \frac{1}{k+1}\max_{\|v\| = 1} \frac{\|(\tilde R^T \tilde R +\lambda I)^{-\frac{1}{2}}(\tilde R^T \tilde R +\lambda I)^{-\frac{1}{2}}v\|^2}{\|(\tilde R^T \tilde R + \lambda I)^{-\frac{1}{2}}v\|^4}, \\
	& \leq & \frac{1}{k+1} \|(\tilde R^T \tilde R +\lambda I)^{-\frac{1}{2}}\|^2\max_{\|v\| = 1} \frac{1}{\|(\tilde R^T \tilde R + \lambda I)^{-\frac{1}{2}}v\|^2},\\
	& = & \frac{1}{k+1} \|(\tilde R^T \tilde R +\lambda I)^{-\frac{1}{2}}\|^2 \|(\tilde R^T \tilde R + \lambda I)^{\frac{1}{2}}\|^2.
\EEAS
The norm of the coefficients $\tilde c_{\lambda}^*$ are thus bounded by
\[
	\|\tilde c^*_\lambda \| \leq \sqrt{\frac{1}{k+1} \frac{\|\tilde R^T \tilde R\| + \lambda}{\lambda}} = \sqrt{\frac{\|\tilde R\|^2 + \lambda}{(k+1)\lambda}}.
\]
We will now bound $\|\Delta \tilde c_\lambda \|$. With equation \eqref{eq:sol_deltactildelambda},
\BEAS
	\|\Delta \tilde c_\lambda\| & = & \left\| \left(I-\frac{(\tilde R^T \tilde R +\lambda I)^{-1}\ones \ones^T}{\ones^T (\tilde R^T \tilde R +\lambda I)^{-1} \ones}\right)(\tilde R^T \tilde R +\lambda I)^{-1}Pc^*_\lambda\right\|, \\
	& \leq &   \left\|I-\frac{(\tilde R^T \tilde R +\lambda I)^{-1}\ones \ones^T}{\ones^T (\tilde R^T \tilde R +\lambda I)^{-1} \ones}\right\|\left\| (\tilde R^T \tilde R +\lambda I)^{-1}\right\|\left\|P \right\|\left\|c^*_\lambda\right\|, \\
	& \leq & \left\|(\tilde R^T \tilde R +\lambda I)^{-1}\right\| \left\|P \right\|\left\|c^*_\lambda\right\|,
\EEAS
where the last inequality is obtained by bounding the norm of a projector. Since $\tilde R^T \tilde R \succeq 0$, we have $ (\tilde R^T \tilde R +\lambda I) \succeq \lambda I $, we get
\[
	\|\Delta \tilde c_\lambda\|  \leq \frac{\|P\|}{\lambda}\|c^*_\lambda\|,
\]
which is the desired result.
\end{proof}

Regularization allows a better control of the impact of perturbations, but also changes the solution $c^*$ into~$c^*_\lambda$. The next section analyses the impact of regularization on the extrapolated solution when there are no perturbations.

\subsection{Regularized Chebyshev Polynomial} \label{ssec:regularized_cheby}
The previous section shows that regularization is important for the control of the perturbations present in \eqref{eq:perturb_LFPI}. However, the convergence analysis becomes more complicated in the perturbation-free case, and we  introduce \emph{regularized Chebyshev polynomials.}
\begin{definition}
The regularized Chebyshev polynomial $C^*_\sigma(x,k,\alpha)$ of degree $k$, range $\sigma$ and regularization parameter $\alpha$ is defined as the solution of
\BEAS
	C^*_\sigma(x,k,\alpha) = \argmin_{C \in \reals_k[x] \,:\, C(1) = 1} \; \max_{x \in [0,\sigma]} C^2(x) + \alpha \|C\|^2,
\EEAS
where $\|C\|$ corresponds to the $\ell_2$ norm of the coefficients of polynomial $C$. We write the maximum value as
\BEA
	S_\sigma(k,\alpha) \triangleq \sqrt{\max_{x \in [0,\sigma]} (C^*_\sigma(x,k,\alpha))^2 + \alpha \| C^*_\sigma(x,k,\alpha) \|^2}. \label{eq:def_s_sigma_k_lambda}
\EEA
\end{definition}
Using this specific polynomial we can now bound the accuracy of the extrapolated point using the regularized algorithm.

\begin{proposition}\label{prop:bound_using_reg_cheby}
Let $c_\lambda^*$ be the output of Algorithm \ref{algo:acc_fixedpoint_reg} using the sequence $x_i$ generated by \eqref{eq:linear_fpi} (with $\|G\| \leq \sigma < 1$) and the parameter $\lambda>0$. The accuracy of the extrapolation is bounded by
\BEQ
	\left\|\sum_{i=0}^k (c^*_\lambda)_i x_i -x^*\right\| \leq  \|(G-I)^{-1}\| \sqrt{S^2_\sigma\left(k, \frac{\lambda}{\|x_0-x^*\|^2}\right) \|x_0-x^*\|^2 - \lambda \|c_\lambda^*\|^2}. \label{eq:bound_reg_cheby_reg}
\EEQ
\end{proposition}
\begin{proof}
Consider the optimization problem in Algorithm \ref{algo:acc_fixedpoint_reg},
\BEAS
	\min_{c^T\ones = 1} \|Rc\|^2 + \lambda \|c\|^2.
\EEAS
Since $r_i = (G-I)(x_i-x^*) = (G-I)G^i(x_0-x^*)$, if we use the polynomial $p$ with coefficients $c$, the problem becomes
\BEA
	 & & \min_{\{ p \in \reals_k[x] :\, p(1) = 1\}} \left\{\|(G-I)p(G)(x_0-x^*)\|^2 + \lambda  \|p\|^2 \right\}, \label{eq:shortcut_minimization_problem} \\
	 	 & \leq & \|x_0-x^*\|^2\min_{\{p \in \reals_k[x] :\, p(1) = 1\}} \left\{\|(G-I)p(G)\|^2 + \frac{\lambda}{\|x_0-x^*\|^2} \|p\|^2\right\},  \nonumber \\
	 	 & \leq & \|x_0-x^*\|^2\min_{\{p \in \reals_k[x]  :\, p(1) = 1\}} \left\{\|G-I\|^2 \|p(G)\|^2 + \frac{\lambda}{\|x_0-x^*\|^2} \|p\|^2\right\},  \nonumber\\
	 	 & \leq & \|x_0-x^*\|^2\min_{\{p \in \reals_k[x]  :\, p(1) = 1\}} \left\{\|p(G)\|^2 + \frac{\lambda}{\|x_0-x^*\|^2} \|p\|^2\right\},  \nonumber
\EEA
where $\|p\|$ is the $\ell_2$ norm of the coefficients of $p$. For simplicity, we write $\bar \lambda = \lambda / \|x_0-x^*\|$ to be the normalized value of $\lambda$. In the optimization problem, since $\|G\| \leq \sigma$, we can consider the worst-case over all symmetric matrices $M$ with $\|M\|\leq \|G\|$ and $M\succeq 0$, written
\BEAS
	\min_{\{ p \in \reals_k[x] :\, p(1) = 1\}}\left\{ \|p(G)\|^2 + \bar\lambda \|p\|^2\right\} \leq
	\min_{\{ p \in \reals_k[x] :\, p(1) = 1\}}\;\max_{M \succeq 0,\; \|M\| \leq \sigma} \left\{\|p(M)\|^2 + \bar \lambda \|p\|^2\right\}.
\EEAS
Because $M$ is symmetric, we only need to look at its eigenvalues which are inside the segment $[0,\sigma]$,
\BEAS
	\min_{\{ p \in \reals_k[x] :\, p(1) = 1\}}\; \max_{M \succeq 0,\; \|M\| \leq \sigma} \left\{\|p(M)\|^2 + \bar \lambda \|p\|^2\right\} & = & \min_{\{ p \in \reals_k[x] :\, p(1) = 1\}} \; \max_{x \in[0,\sigma]}  \left\{p^2(x) + \bar \lambda \|p\|^2\right\}, \\
	& = & S^2_{\sigma}(k,\bar \lambda).
\EEAS
This means that \eqref{eq:shortcut_minimization_problem} is bounded by
\BEA
	\min_{\{ p \in \reals_k[x] :\, p(1) = 1\}} \left\{ \|(G-I)p(G)(x_0-x^*)\|^2 + \lambda  \|p\|^2\right\} \leq \|x_0-x^*\|^2 S^2_{\sigma}(k,\bar \lambda). \label{eq:bound_min_pert_residual}
\EEA
It remains to link the optimization problem to the accuracy of the extrapolation. Indeed,
\BEAS
	\left\|\sum_{i=0}^k (c^*_\lambda)_i x_i -x^*\right\|^2 & = & 	\left\|(G-I)^{-1}\sum_{i=0}^k (c^*_\lambda)_i r_i\right\|^2, \\
	& \leq & \left\|(G-I)^{-1}\right\|^2 \left\|\sum_{i=0}^k (c^*_\lambda)_i r_i\right\|^2, \\
	& = & \left\|(G-I)^{-1}\right\|^2 \left(\left\|\sum_{i=0}^k (c^*_\lambda)_i r_i\right\|^2 + (\lambda - \lambda) \|c^*_\lambda\|^2\right).
\EEAS
By definition, of $c_\lambda^*$,
\BEAS
	\left\|\sum_{i=0}^k (c^*_\lambda)_i r_i\right\|^2 + \lambda \|c^*_\lambda\|^2 & = & \min_{ p \in \reals_k[x]  \,:\, p(1) = 1} \left\{ \|(G-I)p(G)(x_0-x^*)\|^2 + \lambda  \|p\|^2\right\}.
\EEAS
We proved in \eqref{eq:bound_min_pert_residual} that this quantity can be bounded by $S^2_\sigma(k,\bar \lambda) \|x_0-x^*\|^2$, so we finally have
\[
		\left\|\sum_{i=0}^k (c^*_\lambda)_i x_i -x^*\right\| \leq \|(G-I)^{-1}\| \sqrt{S^2_\sigma(k,\bar \lambda) \|x_0-x^*\|^2 - \lambda \|c_\lambda^*\|^2},
\]
which is the desired result.
\end{proof}

\begin{figure}[ht]
	\begin{center}
		\begin{tabular}{cc}
			\psfrag{degree}[t]{Degree $k$}
			\psfrag{lambda}[b]{Regularization parameter $\alpha$}
			\includegraphics[width=0.45\textwidth]{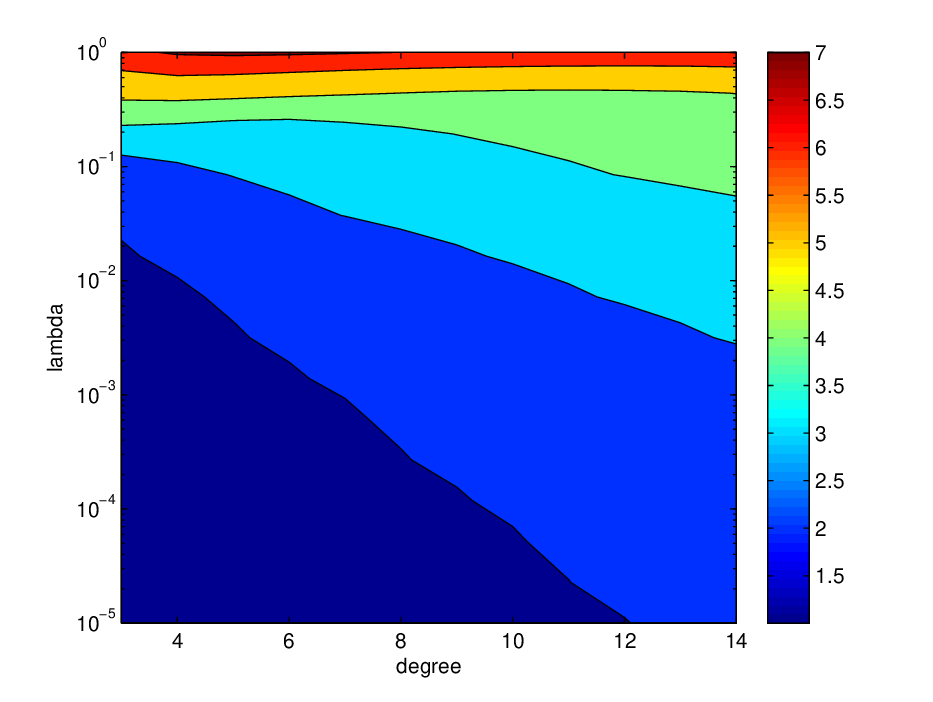}
			&
			\psfrag{degree}[t]{Degree $k$}
			\psfrag{lambda}[b]{Regularization parameter $\alpha$}
			\includegraphics[width=0.45\textwidth]{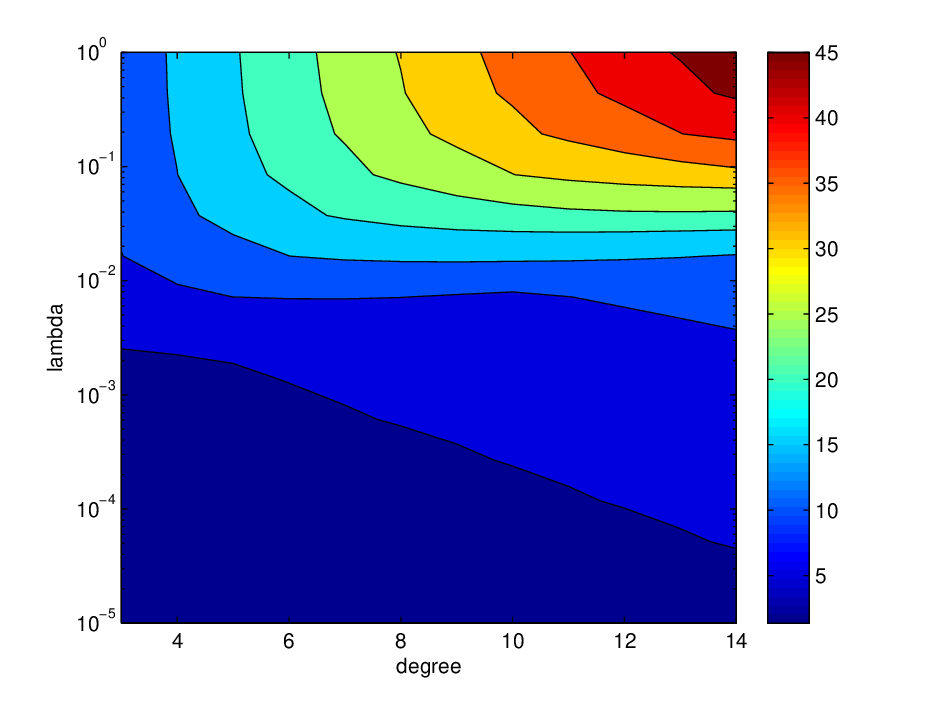}
		\end{tabular}
	\end{center}
	\caption{Ratio between the (worst-case) number of iterations required to reach an arbitrary accuracy using the coefficients of the regularized and non-regularized Chebyshev polynomial for combining the $x_i$. On the left, $\sigma = 0.9$ and on the right $\sigma = 0.999$. We see that the impact of the regularization is more important when $k$ is big, or $\sigma$ close to $1$.\label{fig:comparison_reg_cheby}}
\end{figure}

Regularized Chebyshev polynomials are crucial for the bound on the accuracy of Algorithm \ref{algo:acc_fixedpoint_reg}. Unfortunately, there is no explicit expressions for $S_\sigma(k,\alpha)$ in \eqref{eq:def_s_sigma_k_lambda}. However, this value can be computed numerically using sum-of-squares optimization. We show in Figure \ref{fig:comparison_reg_cheby} the difference of performances when using the coefficients of the regularized Chebyshev polynomial instead of its non-regularized version.

We briefly recall basic results on Sum of Squares (SOS) polynomials and moment problems \citep{Nest00,Lass01,Pari00}, which will allow us to formulate problem~\eqref{eq:def_s_sigma_k_lambda} as a (tractable) semidefinite program. A univariate polynomial is positive if and only if it is a sum of squares. Furthermore, if we let $m(x)=(1,x,\ldots,x^k)^T$
we have, for any $q(x)\in\mathbb{R}_{[2k]}$,
\BEAS
&q(x) \geq 0,~\mbox{for all $x\in\reals$}&\\
&\Updownarrow&\\
&q(x) = m(x)^TCm(x),~\mbox{for some $C\succeq 0$,}&
\EEAS
which means that checking if a polynomial is non-negative on the real line is equivalent to solving a linear matrix inequality (see e.g. \citep[\S4.2]{Bent01} for details). We can thus write the problem of computing the maximum of a polynomial over the real line as
\BEQ\label{eq:sdp-max}
\BA{ll}
\mbox{minimize} & t\\
\mbox{subject to} & t-p(x) = m(x)^TCm(x), \quad\mbox{for all $x\in\reals$}\\
& C\succeq 0,
\EA\EEQ
which is a semidefinite program in the variables $p\in\reals^{k+1}$, $C\in\symm_{k+1}$ and $t\in\reals$, because the first contraint is equivalent to a set of linear equality constraints. Then, showing that $p(x)\geq 0$ on the segment $[0,\sigma]$ is equivalent to showing that the rational fraction
\[
p\left(\frac{\sigma x^2}{1+x^2}\right)
\]
is non-negative on the real line, or equivalently, that the polynomial
\[
(1+x^2)^{k}~p\left(\frac{\sigma x^2}{1+x^2}\right) 
\]
is non-negative on the real line. Overall, this implies that problem~\eqref{eq:def_s_sigma_k_lambda} can be written
\BEQ\label{eq:cheb-reg-sdp}
\BA{rll}
S_{\sigma}(k,\alpha) = & \mbox{min.} & t^2 + \alpha^2\|q\|_2^2\\
&\mbox{s.t.} & t-(1+x^2)^{k+1}~\left(\left(1-\frac{\sigma x^2}{1+x^2}\right)\,q\left(\frac{\sigma x^2}{1+x^2}\right)\right) = m(x)^TCm(x), \quad\mbox{for all $x\in\reals$}\\
&& \ones^Tq=1,\,C\succeq 0,
\EA\EEQ
which is a semidefinite program in the variables $q\in\reals^{k+1}$, $C\in\symm_{k+2}$ and $t\in\reals$.\\

\subsection{Convergence rate} \label{ssec_rate_conv_regalgo}

We will now prove global accuracy bounds, using the following decomposition of the error term,
\BEQ
\sum_{i=0}^k (\tilde c_\lambda^*)_i \tilde x_i - x^* = \underbrace{\sum_{i=0}^k (c_\lambda^*)_i x_i - x^*}_{\textbf{Linear case}} ~+~ \underbrace{\sum_{i=0}^k (\Delta \tilde c_\lambda^*)_i x_i}_{\textbf{Stability}}  ~+~ \underbrace{\sum_{i=0}^k (\tilde c_\lambda^*)_i (\tilde x_i - x_i )}_{\textbf{Nonlinearity}}, \label{eq:decomposition_error}
\EEQ
where $\Delta \tilde c_\lambda^* =\tilde c^*_\lambda - c_\lambda^*$. In the equation above,  the first term is the accuracy of the accelerated method in the noiseless case. The second term corresponds the stability of the coefficients computed by the regularized algorithm when we have some perturbations in the sequence. The last term is the induced error by the combination of the perturbations. The following Theorem shows how to bound these three terms by putting together the results of Propositions \ref{prop:bound_using_reg_cheby} and \ref{prop:stability_reg_algo}.

\begin{theorem}\label{theo:global_bound_acc}
Let $\bar x$ be an arbitrary point in $\mathbb{R}^n$. Given iterates $\tilde x_i$, $i=0, \ldots, k+1$ generated by \eqref{eq:perturb_LFPI}, Algorithm~\eqref{algo:acc_fixedpoint_reg} outputs $x_{\text{extr}} = \sum_{i=0}^k ( \tilde c_\lambda^* )_i \tilde x_i$. Consider the matrices $\bar X$ and $\mathcal{E}$, with columns $\bar X_i = x_i-\bar x$ and $\mathcal{E}_i = \tilde x_i-x_i$ respectively. We have the following bound on the extrapolated point,
\[
	\|x_{\text{extr}} - x^*\| \leq \|x_0-x^*\| S_\sigma(k,\bar \lambda) \sqrt{\kappa^2 + \frac{ \|\bar X\|^2 \|P\|^2}{\lambda^3}} + \frac{\|\mathcal{E}\|}{\sqrt{k+1}} \sqrt{1+\frac{\|\tilde R\|^2}{\lambda}}.
\]
where $\kappa > 1$ with $\|(G-I)^{-1}\|\leq \frac{1}{1-\sigma}= \kappa$.
\end{theorem}
\begin{proof}
The proof is divided into four parts, where the three first parts bound each term of \eqref{eq:decomposition_error} and the last one combines everything. The bound on the first term comes explicitly from Proposition \ref{prop:bound_using_reg_cheby}, with
\BEA
	 \left\|\sum_{i=0}^k (c^*_\lambda)_i x_i -x^*\right\| \leq \kappa \sqrt{S^2_\sigma(k,\bar \lambda) \|x_0-x^*\|^2 - \lambda \|c_\lambda^*\|^2} \label{eq:result_part_one},
\EEA
where $\bar \lambda = \lambda/\|x_0-x^*\|^2$. The second term can be bounded using the fact that both $c_\lambda^*$ and $\tilde c_\lambda^*$ sum to one, so $\Delta \tilde c_\lambda^*$ sum to zero. In this case,
\BEAS
	\left\|\sum_{i=0}^k (\Delta \tilde c_\lambda^*)_i x_i\right\| & = & \left\|\sum_{i=0}^k (\Delta \tilde c_\lambda^*)_i (x_i-\bar x)\right\|, \\
	& \leq & \|\Delta \tilde c_\lambda^*\| \|\bar X\|.
\EEAS
Proposition \ref{prop:stability_reg_algo} bounds the value of $\|\Delta \tilde c_\lambda^*\| $ and yields
\BEA
	\left\|\sum_{i=0}^k (\Delta \tilde c_\lambda^*)_i x_i\right\| \leq \|X\| \frac{\|P\|}{\lambda}\|c_\lambda^*\|. \label{eq:result_part_two}
\EEA
For the third term in~\eqref{eq:decomposition_error}, we have
\BEAS
	 \left\| \sum_{i=0}^k (\tilde c_\lambda^*)_i (\tilde x_i - x_i ) \right\| \leq \| \tilde c_\lambda^*\| \|\mathcal{E}\|.
\EEAS
The norm $\| \tilde c_\lambda^*\|$ can be bounded using Proposition \ref{prop:stability_reg_algo}, with
\BEA
	 \left\| \sum_{i=0}^k (\tilde c_\lambda^*)_i (\tilde x_i - x_i ) \right\| \leq \frac{\|\mathcal{E}\|}{\sqrt{k+1}} \sqrt{1+\frac{\|\tilde R\|^2}{\lambda}}. \label{eq:result_part_three}
\EEA
We finally combine the bounds \eqref{eq:result_part_one}, \eqref{eq:result_part_two} and \eqref{eq:result_part_three} according to the decomposition \eqref{eq:decomposition_error}, to get
\BEA
	& & \left\| \sum_{i=0}^k (\tilde c_\lambda^*)_i x_i-x^* \right\| \leq \label{eq:big_bound_accuracy} \\ 
	& & \kappa \sqrt{S^2_\sigma(k,\bar \lambda) \|x_0-x^*\|^2 - \lambda \|c_\lambda^*\|^2} + \|\bar X\| \frac{\|P\|}{\lambda}\|c_\lambda^*\| + \frac{\|\mathcal{E}\|}{\sqrt{k+1}} \sqrt{1+\frac{\|\tilde R\|^2}{\lambda}}. \nonumber 
\EEA
Here, $\|c_\lambda^*\|$ appears twice in the expression. We remove it by maximizing the bound over $\|c_\lambda^*\|$. The first two terms of~\eqref{eq:big_bound_accuracy} can be written
\[
	x \mapsto \sqrt{a - \lambda x^2} + bx.
\]
with $a = S^2_\sigma(k,\bar \lambda) \|x_0-x^*\|^2 $ and $b=\|\bar X\| \frac{\|P\|}{\lambda} $.  By Proposition \ref{prop:opt_val_sqrt_fun} (in the Appendix), its maximum value is equal to
\[
	\sqrt{a} \sqrt{\kappa^2 + \frac{b^2}{\lambda}},
\]
which is
\[	
	 S_\sigma(k,\bar \lambda) \|x_0-x^*\| \sqrt{\kappa^2 + \frac{ \|\bar X\|^2 \|P\|^2}{\lambda^3}}.
\]
The bound on extrapolation accuracy in~\eqref{eq:big_bound_accuracy} now becomes
\BEAS
	\left\| \sum_{i=0}^k (\tilde c_\lambda^*)_i x_i-x^* \right\| \leq S_\sigma(k,\bar \lambda) \|x_0-x^*\| \sqrt{\kappa^2 + \frac{ \|\bar X\|^2 \|P\|^2}{\lambda^3}} + \frac{\|\mathcal{E}\|}{\sqrt{k+1}} \sqrt{1+\frac{\|\tilde R\|^2}{\lambda}}. 
\EEAS
which is the desired result.
\end{proof}

We can further simplify the bound above by bounding $\|P\|$ using $\sigma$, $\mathcal{E}$ and $\|\bar X\|$.

\begin{proposition}\label{prop:bound_p_r}
	Let $P$ the perturbation matrix defined in \eqref{eq:def_P}. Then
	\BEAS
		\|P\| & \leq & 4(\|\mathcal{E}\|\|R\| + \|\mathcal{E}\|^2), \\
		\|R\| & \leq & \frac{1-\sigma^{k+1}}{1-\sigma}\|x_0-x^*\|,
	\EEAS
	where $\mathcal{E}$ is defined in Theorem \ref{theo:global_bound_acc}, $R$ is the matrix of residuals for the sequence $x_i$ generated by \eqref{eq:linear_fpi}, for $\|G\| \leq \sigma$.
\end{proposition}
\begin{proof}
We begin by the bound on $R$,
\[
	\|R\| \leq \sum_{i=0}^k \|r_i\| \leq \sum_{i=0}^k \|G\|^i\|r_0\| \leq \sum_{i=0}^k \sigma^i\|r_0\|.
\]
Since $r_0 = (G-I)(x_0-x^*)$, and $\|G-I\| \leq 1$, we have $\|r_0\| \leq \|x_0-x^*\|$. Injecting this result in the previous bound gives the desired result,
\[
	\|R\| \leq \sum_{i=0}^k \sigma^i\|x_0-x^*\| = \frac{1-\sigma^{k+1}}{1-\sigma}\|x_0-x^*\|.
\]
Now we prove the bound on $\|P\|$. Let $\tilde R = R + \Delta$ for some perturbation matrix $\Delta$. Then
\BEAS
 	\|P\| & = & \|R^TR-\tilde R^T\tilde R\|, \\
	& \leq & 2\|\Delta\| \| R \| + \|\Delta\|^2.
\EEAS
It remains to bound $\|\Delta\|$. Consider $\bar X$, where each column of $\bar X = (x_i-\bar x)$ for some point $\bar x$. Then we can build $R$ from $X$,
\[
	R = \bar X \begin{bmatrix}
	-1 & 1 \\
	   & -1 & 1 & \\
	   &    & \ddots
	\end{bmatrix} = \bar X D.
\]
It is possible to show $\|D\|\leq 2$. Using the same logic, we can build $\tilde R$,
\[
	\tilde R = (X + \mathcal{E}) D = R + \mathcal{E} D.
\]
By identification, we have $\Delta = \mathcal{E} D$, so $\|\Delta \| \leq \|D\| \|\mathcal{E}\| \leq 2\|\mathcal{E}\|$.
\end{proof}

Assuming again $\|G\| \leq \sigma$, the following propositions bound $\|\bar X\|$ when $\bar x = x^*$.

\begin{proposition} \label{prop:bound_barx} Let $\bar X$ be the matrix built with the columns $\bar X_i = x_i-\bar x$, where the sequence $x_i$ is generated by \eqref{eq:linear_fpi} and $\bar x = x^*$. If $ \|G\| \leq \sigma$, where $G$ is the matrix present in \eqref{eq:linear_fpi}, the norm of $\bar X$ is bounded by
\BEA
	\|\bar X\| \leq \frac{1-\sigma^{k+1}}{1-\sigma}\|x_0-x^*\|.
\EEA
\end{proposition}
\begin{proof}
Since each column of $\bar X$ correspond to $x_i-x^*$,
\BEAS
	\|\bar X\| & \leq & \sum_{i=0}^k \|x_i-x^*\| \leq \sum_{i=0}^k \|G^i(x_0-x^*)\| \leq \sum_{i=0}^k \|G\|^i \|x_0-x^*)\|.
\EEAS
Because $\|G\| \leq \sigma < 1$,
\[
	\|\bar X\| \leq \frac{1-\sigma^{k+1}}{1-\sigma} \|x_0-x^*\|,
\]
which is the desired result.
\end{proof}

The bound of Theorem \ref{theo:global_bound_acc} is quite generic. For now, we only need a sequence $\tilde x_k$ generated by a perturbed fixed-point process which is convergent and differentiable, and the accuracy depends on matrices $\tilde R$ and $\mathcal{E}$. The next section will bound these quantities when the fixed point process is the gradient descent algorithm.

\subsection{Accelerating Gradient Descent}
Assume the sequence $\tilde x_i$ is generated by the gradient descend algorithm,
\[
	\tilde x_{i+1} = \tilde x_i -\frac{1}{L}f'(\tilde x_i),
\]
where $f$ is a $\mu$-strongly convex, $L$-smooth function with a Lipschitz-continuous Hessian with constant $M$. In this case, we can bound the values $\|\tilde R\|$ and $\|\mathcal{E}\|$ and hence $\|P\|$. We show the following result in Section \ref{app:proof_explicit_bounds_gradient}.

\begin{proposition} \label{prop:explicit_bounds_gradient}
When using gradient method on a $\mu$-strongly convex, $L$-smooth function with a Lipschitz-continuous Hessian with constant $M$, we have the following bounds,
\BEA
	\|\tilde R\| &\leq & \frac{1-\sigma^{k+1}}{1-\sigma} L\|x_0-x^*\|,\\
	\|\mathcal{E}\| &\leq & (k+2)^2 \frac{M}{4L}\|x_0-x^*\|^2,
\EEA
where $\sigma = 1-\frac{\mu}{L}$ satisfies $\|G\| \leq \sigma$.
\end{proposition}

\begin{figure}
	\psfrag{degree}[t][b]{Degree $k$}
	\psfrag{speedup}[b][t]{Speedup}
	\psfrag{Nesterov-reference}{\footnotesize Nesterov}
	\psfrag{normalized-rate-extrapolation}{\footnotesize Regularized acceleration}
	\psfrag{Gradient method}{\footnotesize Gradient method}
	\includegraphics[width=0.50\textwidth]{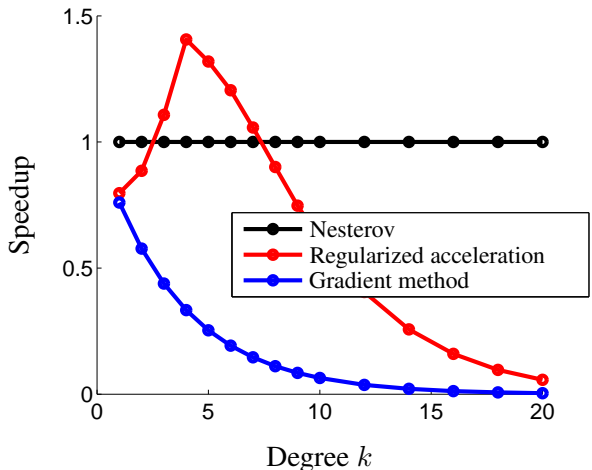}
	\caption{Convergence speedup relative to Nesterov's accelerated method of theoretical bound in Theorem~\ref{theo:global_bound_acc} and gradient method, using upper bounds from Propositions \ref{prop:bound_p_r}, \ref{prop:bound_barx} and \ref{prop:explicit_bounds_gradient}. We see that our (highly conservative) bound shows a slight speedup when $k$ is well chosen.}
	\label{fig:speedup_acc}
\end{figure}

Using these expressions, we can compare convergence rates between convergence acceleration in~Algorithm~\ref{algo:acc_fixedpoint_reg} and Nesterov's method. In Figure \ref{fig:speedup_acc} we illustrate the difference on a particular instance where $\|x_0-x^*\|=10^{-4}$, $L = 1$, $\mu = M = 0.1$. We see that, despite the highly conservative nature of this bound, for small $k$ at least our method is faster than Nesterov's acceleration.

When using the gradient method, this result bounds all quantities present in Theorem \ref{theo:global_bound_acc} as a function of $\mu,\,L,\,M$ and $\|x_0-x^*\|$.  Asymptotically, i.e. when  ${\|x_0-x^*\| \rightarrow 0}$ and we are starting close enough to the optimal point, we show that we recover the acceleration rate for linear sequences in Proposition~\ref{prop:linear_rate_acceleration} if the regularization parameter $\lambda$ il well-chosen.

\begin{proposition}
Assume we used the gradient method on a $L$-smooth and $\mu$-strongly convex function with Lipschitz-continuous Hessian to generate the sequence $\tilde x_i$. Setting $\lambda =  O(\|x_0-x^*\|^s)$ with $s\in ]2,\frac{8}{3}[$ and ${\|x_0-x^*\| \rightarrow 0}$ (i.e., starting close enough to the optimal point), then the rate of convergence of the extrapolated point is bounded by
\[
	\lim\limits_{\|x_0-x^*\| \rightarrow 0} \frac{\|\sum_{i=0}^k \tilde c_{\lambda}^* \tilde x_i - x^*\|}{\|x_0-x^*\|} \leq \kappa \dfrac{2\beta^k}{1+\beta^{2k}} \,, \qquad \beta = \frac{1-\sqrt{\kappa^{-1}}}{1+\sqrt{\kappa^{-1}}},
\]
where $\kappa = \frac{L}{\mu}$.
\end{proposition}
\begin{proof}
The bounds above show
\[
	\|\bar X\| = \|\tilde R\| = O(\|x_0-x^*\|) \;, \qquad \|\mathcal{E}\| = O(\|x-x^*\|^2) \;, \qquad \|P\| = O(\|x_0-x^*\|^3) .
\]
Let $\lambda = O(\|x_0-x^*\|^s)$ for some scalar $s$. The bound of Theorem \ref{theo:global_bound_acc} normalized by $\|x_0-x^*\|$ becomes
\[
	S_{\sigma}(k,O(\|x_0-x^*\|^{s-2})) \sqrt{\kappa^2 + O(\|x_0-x^*\|^{8-3s}) } + \sqrt{O(\|x_0-x^*\|^2) + O(\|x_0-x^*\|^{4-s})}.
\]
If $4-s < 0$, clearly the last terms vanishes when $\|x_0-x^*\|\rightarrow 0$. It remains to analyze
\[
	\lim\limits_{\|x_0-x^*\| \rightarrow 0} \; S_{\sigma}(k,O(\|x_0-x^*\|^{s-2})) \sqrt{\kappa^2 + O(\|x_0-x^*\|^{8-3s}) }.
\]
If $s \in ]2,8/3[$ then $s-2 > 0$ and $8-3s > 0$, implying $(\|x_0-x^*\|^{s-2}) \rightarrow 0$ and $O(\|x_0-x^*\|^{8-3s}) \rightarrow 0$ when $\|x_0-x^*\|\rightarrow 0$. The bound finally becomes
\[
	\lim\limits_{\|x_0-x^*\| \rightarrow 0} \frac{\|\sum_{i=0}^k \tilde c_{\lambda}^* \tilde x_i - x^*\|}{\|x_0-x^*\|} \leq \kappa S_{\sigma}(k,0).
\]
However, $ S_{\sigma}(k,0) $ is exactly equal to the maximum value of the rescaled (non-regularized) Chebyshev polynomial $T_{k}(x,\sigma)$, so by equation \eqref{eq:chebysol},
\[
	S_{\sigma}(k,0) = \max_{x\in [0,\sigma]}T_{k}(x,\sigma) =  \dfrac{2\beta^k}{1+\beta^{2k}}.
\]
This result conclude the proof.
\end{proof}

In other words, the result above means that the bound of Theorem~\ref{theo:global_bound_acc} tends to be the bound of Proposition~\ref{prop:linear_rate_acceleration} (for the case where $k<m$) and, asymptotically, we recover the optimal rate of convergence  in \citep{Nest03a}. In fact, the result may also hold for other kinds of methods, because the proof only needs
\[
	\|\bar X\| = \|\tilde R\| = O(\|x_0-x^*\|) \;, \qquad \|\mathcal{E}\| = O(\|x-x^*\|^2) \;, \qquad \|P\| = O(\|x_0-x^*\|^3) .
\]
These assumptions are not too restrictive, and are often encountered when using deterministic, twice-differentiable, linearly convergent iterations $g$ in \eqref{eq:iter}. 

In practice of course, $\|x-x^*\|$ is unknown and the regularization parameter $\lambda$ should decrease fast enough to ensure $S_\sigma(k,\bar \lambda) \rightarrow S_\sigma(k,0)$, but not too fast otherwise the algorithm becomes unstable. An adaptive strategy thus ensures a good convergence rate, which is what we detail next.

\subsection{Adaptive regularization}
The major problem of the regularized Algorithm \ref{algo:acc_fixedpoint_reg} is the presence of the parameter $\lambda$, unknown in advance. Of course, one can use the bound in Theorem \ref{theo:global_bound_acc} to search the best~$\lambda$, but this requires a lot of information on the problem, like the constants $L,\,\mu$ and $M$ as well as the distance to the optimum $\|x_0-x^*\|$. Moreover, the bound is extremely pessimistic and does not correspond to the good numerical performances of the algorithm.

To avoid this problem we use adaptive strategy to find $\lambda$, based on grid search, which requires $k$ additional calls to $f(x)$. In comparison, we also need to call $k$ times the oracle for common adaptive strategy in the (accelerated) gradient method. For example, the backtracking line-search over the constant $L$ requires the evaluation of $f(x_i)$ at each iteration $i=1...k$.

Finally, the introduction of the regularization parameter introduces some damping in the acceleration algorithm, in the sense that the step length $x_{\text{extr}}(\lambda) - x_0$ is reduced with higher values of $\lambda$. A simple line search over the step-size, which consists in finding a good scalar $t$ which minimizes the function, solving
\[
	\min_{t>0} f\big( x_0 + t ( \underbrace{ x_{\text{extr}}(\lambda) - x_0 }_{\text{Extrapolation step}} ) \big)
\]
significantly improves the solution. Nevertheless, this requires further calls to $f(x)$, and an inexact line-search is usually preferable. We start with $t=1$, then multiply the value by two until the objective function increases,
\[
	f\big( x_0 + t (  x_{\text{extr}}(\lambda) - x_0 ) \big) < f\big( x_0 + 2t (  x_{\text{extr}}(\lambda) - x_0 ) \big).
\]
In our numerical experiments, this line-search dramatically increases acceleration performances.

We summarize all the steps detailed above as the \emph{Adaptive Regularized Convergence Acceleration} Algorithm~\ref{algo:adaptive_blackbox_acceleration}. The only required inputs are the sequence $\tilde x_i$ generated by the optimization algorithm and the objective function $f$. 

\begin{algorithm}[h]
\caption{Adaptive Regularized Nonlinear Acceleration of Convergence}
\label{algo:adaptive_blackbox_acceleration}
\begin{algorithmic}[1]
\REQUIRE Sequence $\{\tilde x_0, \tilde x_1,..., \tilde x_{k+1}\}$, bounds $[\lambda_{\min},\lambda_{\max}]$, objective function $f(x)$.
\STATE Divide the segment $[\lambda_{\min},\lambda_{\max}]$ into $k$ points $\{\lambda_j\}$ using a logarithmic scale.
\STATE Compute the residual matrix $\tilde R$ such that $\tilde R_i = \tilde x_{i+1} - \tilde x_{i}$.
\STATE Build the matrix $M = \tilde R^T \tilde R / \|\tilde R^T \tilde R\|$
\FOR{$j$ in $1...k$}
\STATE Solve in $z$ the linear system $(M+\lambda_j)z = \ones $
\STATE Normalize the solution, $\tilde c_{\lambda_j}^* = {z}/{\ones^Tz}$
\STATE Compute $x_{\text{extr}}(\lambda_j) = \sum_{i=0}^k (\tilde c_{\lambda_j}^*)_{_i} \tilde x_i$
\ENDFOR
\STATE Pick $x_{\text{extr}}^* = \argmin_{j=1..k} f(x_{\text{extr}}(\lambda_j))$
\STATE Define $F_t = f(x_0 + t(x_{\text{extr}}^* - x_0))$
\STATE Initialize with $t=1$
\WHILE{$F_{2t} < F_t$}
\STATE Update $ t = 2t $
\ENDWHILE
\ENSURE Return $(x_0 + t(x_{\text{extr}}^* - x_0))$, the extrapolated point.
\end{algorithmic}
\end{algorithm}

\subsection{Computational Complexity of Convergence Acceleration}\label{ss:complex}
In Algorithm \ref{algo:acc_fixedpoint_reg}, computing the coefficients $\tilde c^*_\lambda$ means solving the $k\times k$ system $(\tilde R^T \tilde R + \lambda I)z = \ones$. We then get $\tilde c^*_\lambda = z/(\ones^T z)$. This can be done in both batch and online mode. We will see that, in any case, we end with a complexity of $O(nk^2 + k^3)$, for a small value of $k$ (usually, $k=5$). The complexity of the acceleration algorithm is linear in the dimension, thus adding a negligible additional computation cost to the original procedure.

\subsubsection{Online updates.}
Here, we receive the vectors $r_i$ one by one from the optimization algorithm, and we would like to solve the linear system in parallel of the optimization algorithm. In this case, we perform low-rank updates on the Cholesky factorization of the system matrix. At iteration $i$, we have the Cholesky factorization $LL^T = \tilde R^T \tilde R + \lambda I$, where $L$ is a triangular matrix. We receive a new vector $r_+$ and we want 
\[ 
L_+L^T_+ = \begin{bmatrix}
L & 0 \\
a^T & b
\end{bmatrix}\begin{bmatrix}
L^T & a \\
0 & b
\end{bmatrix} = \begin{bmatrix}
\tilde R^T \tilde U + \lambda I & \tilde R^T r_+\\
(\tilde R^Tr_+)^T & r_+^Tr_+ + \lambda
\end{bmatrix} .
\]
We can explicitly solve this system in variables $a$ and $b$, and the solutions are
\[
a = L^{-1}\tilde R^Tr_+  ,\quad  b = a^Ta+\lambda.
\]
The complexity of this update is thus $O(i\,n + i^2)$, i.e. the matrix-vector multiplication of $\tilde R^Tr_{+}$ with cost $O(i\,n)$ and solving a $i\times i$ triangular system with cost $O(i^2)$. Since we need to do it $k$ times, the final complexity is thus $O(nk^2 + k^3)$. 

\subsubsection{Batch mode.}
The complexity is divided in two parts: First, we need to build the linear system itself. Since $\tilde R\in\mathbb{R}^{n\times k}$, it takes $O(nk^2)$ flops to perform the multiplication. Then we need to solve the linear system $(\tilde R^T\tilde R + \lambda I)z = \ones$ which can be done by Gaussian elimination (in particular when $k$ is small), by Cholesky factorization or by using an iterative method like conjugate gradient. It takes $O(k^3)$ flops to solve the linear system in the worst case, meaning that the overall complexity is $O(nk^2+k^3)$.

\section{Extensions \& Links with other Methods}\label{sec:ext}

\subsection{Smooth Minimization}\label{ss:smooth}
We can extend our results to smooth functions that are not strongly convex using a simple regularization trick which we trace back at least to \citep{Haza14}. Suppose we seek to solve 
\[
\min_{x\in\reals^n} f(x)
\]
in the variable $x\in\reals^n$, where $f(x)$ has a Lipschitz continuous gradient with parameter $L$ with respect to the Euclidean norm, but is not strongly convex. Assume for simplicity that the initial iterate $x_0$ is close enough to the optimum so that $D\triangleq \|x_0-x^*\| \geq \|x_k-x^*\|$ for any $k\geq 0$. We can approximate the above problem by 
\BEQ\label{eq:approx-strong}
\min_{x\in\reals^n} f_\varepsilon(x) \triangleq f(x) + \frac{\varepsilon}{2D^2} \|x_0-x\|_2^2
\EEQ
in the variable $x\in\reals^n$, where $f_\varepsilon(x)$ has a Lipschitz continuous gradient with parameter $L+\varepsilon/D^2$ with respect to the Euclidean norm, is strongly convex with parameter $\varepsilon/D^2$ with respect to the same norm. Furthermore $f_\varepsilon(x)$ is an $\varepsilon$ approximation of $f(x)$ near the optimum and we get
\BEAS
f(x_k)-f(x^*) & = & f_\varepsilon(x_k) - \frac{\varepsilon}{2D^2} \|x_0-x_k\|_2^2 - f_\varepsilon(x^*) + \frac{\varepsilon}{2D^2} \|x_0-x^*\|_2^2 ,\\
& \leq & f_\varepsilon(x_k) - f_\varepsilon(x^*) + \frac{\varepsilon}{2} ,\\
& \leq & f_\varepsilon(x_k) - f_\varepsilon(x_\varepsilon^*) + \frac{\varepsilon}{2},
\EEAS
using the smoothness of $f_\varepsilon(x)$ and writing $x_\varepsilon^*$ the optimum of problem~\eqref{eq:approx-strong}. It suffices to optimize $f_\varepsilon$ up to $\varepsilon/2$ to find an $\varepsilon$-solution for the original problem. The linear convergence of gradient \citep{Nest03a} algorithms guarantees
\[
	f_\varepsilon(x_k) - f_\varepsilon(x_\varepsilon^*) = \frac{(L+\varepsilon)D^2}{2} r^k\,, \qquad r=1-\frac{2\varepsilon}{LD^2+2\varepsilon}.
\]
The number of iterations required to reach a target precision $\varepsilon/2$ is thus bounded by
\[
k = O\left(\frac{\log((L+\varepsilon)D^2/\varepsilon)}{\log(1/r)}\right).
\]
By replacing the value of $r$, we have
\[
\log(1/r) \sim 1-r = \frac{LD^2}{\varepsilon},
\]
while accelerated algorithms have $r=1-\sqrt{\varepsilon/(LD^2+\varepsilon)}$ which yields
\[
\log(1/r) \sim \sqrt{\frac{LD^2}{\varepsilon}}.
\]
Up to a logarithmic constant, these upper bounds match the complexity of gradient and accelerated gradient methods. Overall, an algorithm for strongly convex function used with this regularization trick recovers an $\varepsilon$-approximated solution. This means we can always reduce a not strongly convex problem to \eqref{eq:fprob}, where our acceleration analysis applies.

\subsection{Convergence Acceleration on Gradient Method for Quadratic Functions}
Assume we want to minimize a quadratic function $f$. Its gradient reads, for $A \in \reals^{n\times n}$ a symmetric positive definite matrix,
\[
	\nabla f (x) = A (x-x^*).
\]
This formulation is equivalent to $\nabla f = Ax-b$, where $b = Ax^*$ but it will be more convenient in this section to manipulate directly $x^*$. Let $\mu I \preceq A \preceq L I$ so that the function $f$ is strongly convex of constant $\mu$ and smooth of constant $L$. If we use the fixed-step gradient method, with step-size $1/L$, 
\BEA\label{eq:gradient_iteration}
	x_{i+1} = x_k - \frac{1}{L} \nabla f(x_k) = x_k - \frac{1}{L} A(x_k-x^*).
\EEA
The fixed point iteration corresponds to
\[
	g(x) = (I-A/L)(x-x^*) + x^*.
\]
Notice that $g(x_{k+1})$ and \eqref{eq:gradient_iteration} are equivalent. The Jacobian of $g$ is thus equal to $(I-A/L)$. We have the following bounds on $G$,
\[
	0 \preceq G \preceq \left(1-\frac{\mu}{L}\right)I.
\]
By consequence, $\sigma = 1-\frac{\mu}{L}$, thus the rate of convergence of our method is linear and the bound is
\[	
	\|x_k-x^*\| \leq \left(1-\frac{\mu}{L}\right)^k\|x_0-x^*\|.
\]
However, if we use Algorithm \eqref{algo:acc_fixedpoint}, we combine the iterates $x_i$ with coefficients $c^*$ (computed by formula \eqref{eq:linearsystem_cstar}). By equations \eqref{eq:rate_linear} and \eqref{eq:error_extrapolation} the accuracy of this extrapolation is bounded by
\BEA
	\left\|\sum_{i=0}^N c_i^*x_{i} - x^*\right\| \leq \frac{L}{\mu} \frac{2\beta^k}{1+\beta^{2k}} \|x_0-x^*\|, \qquad \text{where } \; \beta = \frac{1-\sqrt{\frac{\mu}{L}}}{1+\sqrt{\frac{\mu}{L}}}. \label{eq:rate_acc_gradient}
\EEA
This bound matches the rate obtained using the optimal method in \citep{Nest03a}. Outside of the normalization constraint, this is very similar to the convergence analysis of Lanczos' method.

\subsection{Convergence acceleration versus conjugate gradient} \label{ssec:conjgrad} The rate of convergence obtained above also matches that of the conjugate gradient within a factor $L/\mu$. Indeed, the acceleration algorithm has a strong link with the conjugate gradient. Denote $\|v\|_M= \sqrt{v^TMv}$ the norm induced by the positive definite matrix $M$. Also, assume we want to solve $Ax=b$ using conjugate gradient method (where $A$ is assumed to be symmetric and positive definite). By definition, at the $k$-th iteration, the conjugate gradient computes an approximation of $x^*$ which follows
\[
\mathop{\mathrm{argmin}}\limits_{x \in \mathcal{K}_k} \|x-x^*\|_A,
\]
where $\mathcal{K}_k = \Span\{b,Ab,...,A^{k-1}b\} = \Span\{Ax^*,A^2x^*,...,A^{k}x^*\}$ is called a Krylov subspace. Since the constraint $x\in\mathcal{K}_k$ impose us to build $x$ from a linear combination of the basis of $\mathcal{K}_k$, we can write
\[
x = \sum_{i=0}^{k-1}c_iA^{i+1}x^* = q(A)x^*,
\] 
where $q(x)$ is a polynomial of degree $k$ and $q(1) = 0$. So the  conjugate gradient method solves
\[ 
\mathop{\mathrm{argmin}}\limits_{\{q\in\reals_k[x]:\,q(0)=0\}} \|q(A)x^*-x^*\|_A  = \mathop{\mathrm{argmin}}\limits_{\{\hat q\in\reals_k[x]:\,\hat q(0)=1\}} \|\hat q(A)x^*\|_A,
\]
which is very similar to the equations in \eqref{eq:min_polynomial}. However, while conjugate gradient has access to an oracle giving the result of the product between $A$ and any vector $v$, the acceleration algorithm can only use the iterations produced by \eqref{eq:linear_fpi}, so it does not require the knowledge of $A$. Moreover, the convergence of conjugate gradient is analyzed in another norm ($\|\cdot\|_A$ instead of $\|\cdot\|_2$), which explains why a condition number appears in the bound \eqref{eq:rate_acc_gradient}.

Analysis of convergence on conjugate gradient often use Chebyshev's polynomial, like the acceleration algorithm \eqref{algo:acc_fixedpoint}. We will now see that Nesterov's algorithm generates also a polynomial, making the convergence analysis for quadratics easier.

\subsection{Chebyshev's Acceleration and Nesterov's Accelerated Gradient Method}\label{ss:cheby}
In Proposition \ref{prop:rateconv_linearacc}, we bounded the rate of convergence of Algorithm \ref{algo:acc_fixedpoint} using Chebyshev polynomials. In fact, this is exactly the idea behind Chebyshev's semi-iterative method, which uses these coefficients in order to accelerate gradient descent on quadratic functions. Here, we present Chebyshev semi-iterative acceleration and its analysis, then  use the same arguments on Nesterov's method. These points were also discussed in~\citep{Hard13}.

Assume as above that we use the gradient method to minimize a quadratic function, we get the recurrence \eqref{eq:gradient_iteration}. We see easily that 
\[
	x_k = x^* + G^k(x_0-x^*).
\]
Since $\|G\|_2\leq 1-\frac{\mu}{L} = \sigma$, the rate of convergence is $\|x_k-x^*\|_2 \leq \sigma^k\|x_0-x^*\|_2$. Moreover, if we average the vectors $x_i$ using coefficients $c_i$ (with unitary sum) from 0 to $k$, we get
\[
	\sum_{i=0}^{k} c_ix_i = x^*+p(G)(x_0-x^*)
\]
for $p \in \reals_k[x]$ a polynomial of coefficients $c$. Instead of using Algorithm \eqref{algo:acc_fixedpoint}, which minimizes the combination of the residuals instead of the error term, we will use the coefficients of the rescaled Chebyshev polynomial \eqref{eq:rescaledcheby}. Recall this polynomial makes $\|p(G)\|_2$ small for all matrices $G$ such that $0\preceq G\preceq \sigma I$. In other terms, the rescaled Chebyshev polynomial satisfies
\BEAS
	T(x) & = & \arg\min_{\substack{p\in\reals[x]\\p(1)=1}}~\max_{0\preceq G \preceq \sigma I}\left\|p\left(G\right)\right\|_2,\\
	& = & C_k(t_\sigma(x)).
\EEAS
where $T_k$ and $t_\sigma$ are also defined in \eqref{eq:rescaledcheby}. Furthermore, the Chebyshev polynomials can be constructed using a three-terms recurrence
\[
C_k(x) = xC_{k-1}(x)-C_{k-2}(x).
\]
The same holds for $T_k(x)$, with
\BEAS
\alpha_k &=& t(1)\alpha_{k-1}-\alpha_{k-2}, \\
z_{k-1} &=& y_{k-1} - \nabla f(y_{k-1}),\\
y_{k} &=&  \frac{\alpha_{k-1}}{\alpha_{k}}\left(\frac{2z_{k-1}}{\sigma}-y_{k-1} \right)-\frac{\alpha_{k-2}}{\alpha_{k}}y_{k-2}.
\EEAS
This scheme looks very similar to Nesterov's accelerated gradient method, which reads
\BEAS
z_{k-1} &=& y_{k-1} - \nabla f(y_{k-1})\\
y_{k} &=& z_{k-1} + \beta_k(z_{k-1}-z_{k-2})
\EEAS
Compared with Chebyshev acceleration, Nesterov's scheme is iteratively building a polynomial $N_k(x)$ with $y_k-y^* = N_k\left(G\right)(y_0-x^*)$. If we replace $z_{k}$ by its definition in the expression of $y_k$ in the Nesterov's scheme we get the following recurrence of order two
\BEAS
y_{k}-x^* &=& (1+\beta_k)G(y_{k-1}-y^*) - \beta_k G(y_{k-2}-y^*),\\
&=& G\left((1+\beta_k)N_{k-1}\left(G\right) - \beta_k N_{k-2}\left(G\right)\right)(y_0-x^*).
\EEAS
We can extract the polynomial $N_k$, which reads
\[
N_k(x) = x((1+\beta_k)N_{k-1}(x) - \beta_k N_{k-2}(x)),
\] 
with initial conditions $N_0(x) = 1$ and $N_1(x) = x$. Notice that $N_k(1) = 1$ for all $k$. 

When minimizing smooth strongly convex functions with Nesterov's method, we use 
\[
\beta_k = \frac{\sqrt{L}-\sqrt{\mu}}{\sqrt{L}+\sqrt{\mu}}.
\] 
Moreover, empirically at least, the maximum value of $N_k(x)$ in the interval $[0,\sigma]$ is $N_k(\sigma)$. We conjecture that this always holds. We thus have the following recurrence
\BEAS
N_{k}(\sigma) &=& \sigma\left(\left(1+\beta\right)N_{k-1}(\sigma) - \beta N_{k-2}(\sigma)\right)
\EEAS
To get linear convergence with rate $r$, we need $N_k\leq rN_{k-1} \leq r^2N_{k-2}$, or again
\[
N_{k}(\sigma) \leq \sigma\left(\left(1+\beta\right)rN_{k-2}(\sigma) - \beta N_{k-2}(\sigma)\right) =\sigma\left(\left(1+\beta\right)r - \beta \right) N_{k-2}(\sigma).
\]
Now, consider the condition
\[ 
\sigma\left(\left(1+\beta\right)r - \beta \right) \leq r^2.
\]
We have that Nesterov's coefficients and rate, i.e. $\beta={(1-\sqrt{\mu/L})}/{(1+\sqrt{\mu/L})}$ and $r = (1-\sqrt{\mu/L})$, satisfy this condition, showing that Nesterov's method converges with a rate at least $r=(1-\sqrt{\mu/L})$ on quadratic problems. This provides an alternate proof of Nesterov's acceleration result on these problems using Chebyshev polynomials (provided the conjecture on $N(\sigma)$ holds).

\section{Numerical Experiments}\label{sec:numres}
In this section, we evaluate the performance of the adaptive acceleration methods without/with line-search on the step size, described in Algorithm \ref{algo:adaptive_blackbox_acceleration}.

\subsection{Minimizing logistic regression}
We begin by testing our methods on a regularized logistic regression problem written
\[
f(w) = \sum_{i=1}^m \log\left( 1+\exp(-y_i\xi_i^Tw) \right)+ \frac{\tau}{2}\|w\|^2_2 ,
\]
where $Z=[\xi_1,...,\xi_m]^T\in\reals^{m\times n}$ is the design matrix and $y$ is a $\{-1,1\}^m$ vector of labels. The Lipschitz constant of the logistic regression is $L = \|Z\|_2^2/4+\tau$ and the strong convexity parameter is $\mu = \tau$. We solve this problem using several algorithms.
\begin{itemize}\itemsep +1ex
	\item Fixed-step gradient method for smooth strongly convex functions \citep[Th.\,2.1.15]{Nest03a} 
	\[
	x_{k+1} = x_k-\frac{2}{L+\mu}\nabla f(x_k).
	\]
	\item Accelerated gradient method for smooth strongly convex functions \citep[Th.\,2.2.3]{Nest03a}
	\BEAS
	x_{k+1} &=& y_k-\frac{1}{L}\nabla f(y_k),\\
	y_{k+1} &=& x_{k+1} + \frac{\sqrt{L}-\sqrt{\mu}}{\sqrt{L}+\sqrt{\mu}}\left( x_{k+1}-x_k \right).
	\EEAS
	\item The accelerated gradient method with backtracking line-search on the parameter $L$.
	\item The Adaptive acceleration algorithm \ref{algo:adaptive_blackbox_acceleration} on $k$ iterations of gradient descent without line search (written RNA k).
	\item The Black-box acceleration algorithm \ref{algo:adaptive_blackbox_acceleration} (written RNA k + LS) on $k$ iterations of gradient descent.
\end{itemize} 
The matrix $Z$ is build using datasets \emph{Sonar} (60 features, 208 points), \emph{Madelon} (500 features, 4400 points) or \emph{Sido0} (4932 features, 12678 points), concatenated with a column of ones. The optimization is done on the raw data, i.e. without normalization. The starting point is always $w_0 = 0$.\\

\begin{figure}[h]
\begin{center}
  	\psfrag{valf}[b][t]{$f(x_k)-f(x^*)$}
	\psfrag{xlabel}[t][b]{Gradient oracle calls}
	\psfrag{Grad}{\scriptsize Gradient}
   	\psfrag{GradNest}{\scriptsize Nesterov}
   	\psfrag{RMPE10}{\scriptsize RNA 5}
    \psfrag{AMPE10}{\scriptsize Acc. 5}
    \includegraphics[width=0.5\textwidth]{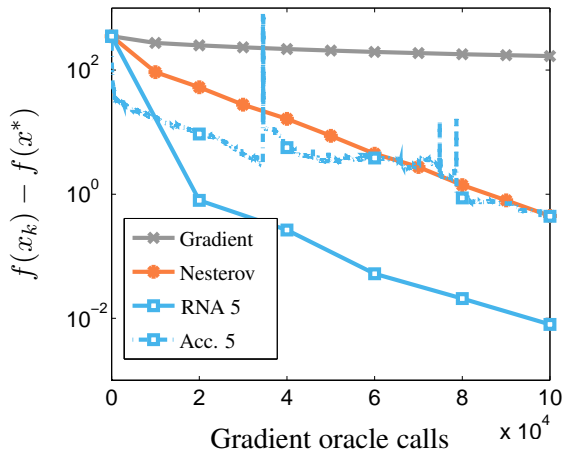}
\end{center}
   \caption{Logistic regression on \emph{Madelon UCI Dataset} with a condition number equal to $1.2\cdot 10^{9}$, solved using Gradient method, Nesterov's method and two versions of the acceleration algorithm applied to the gradient descent: the acceleration algorithm \ref{algo:acc_fixedpoint} (called \textit{Acc. 5})  and the adaptive Regularized Nonlinear Acceleration algorithm \ref{algo:adaptive_blackbox_acceleration} (called \textit{RNA 5}) applied to 5 iterations of the gradient descend. We see that without regularization, the acceleration is unstable because $\|(\tilde R^T \tilde R)^{-1}\|_2$ is huge (cf. Proposition~\ref{prop:sensitivity_acc}).\label{fig:unstableAMPE}}
\end{figure}

\begin{figure}[h]
	\begin{center}
		\begin{tabular}{cc}
		  	\psfrag{valf}[b][t]{$f(x_k)-f(x^*)$}
			\psfrag{xlabel}[t][b]{Gradient oracle calls}
			\psfrag{Grad}{\scriptsize Gradient}
			\psfrag{GradNest}{\scriptsize Nesterov}
			\psfrag{GradNestBack}{\scriptsize Nest. + back.}
			\psfrag{RMPE5}{\scriptsize RNA 5}
			\psfrag{RMPE5LS}{\scriptsize RNA 5 + LS}
			\includegraphics[width=0.45\textwidth]{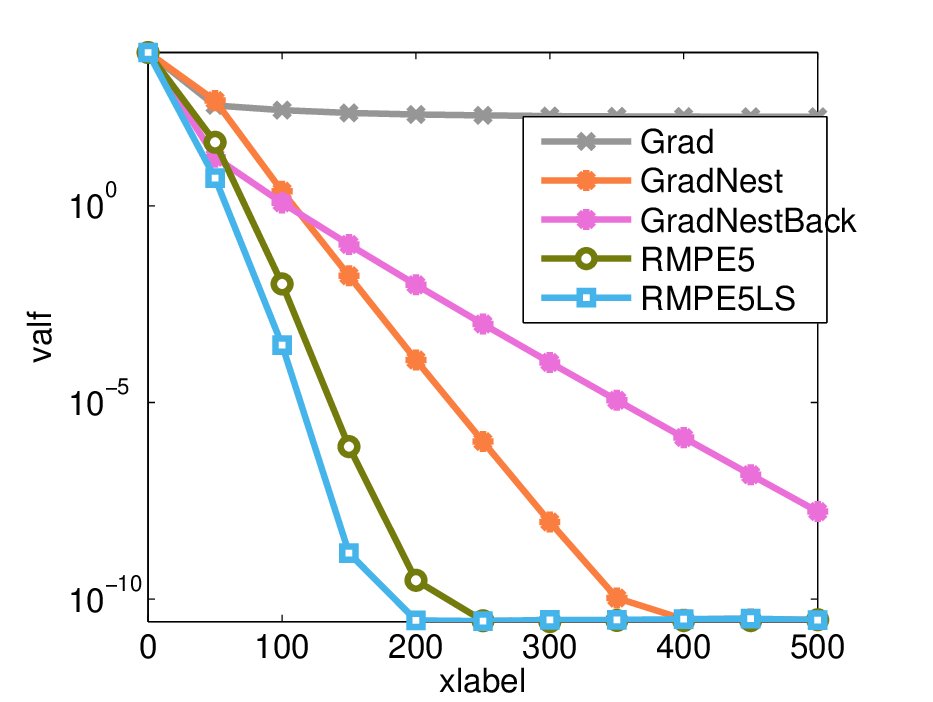}
		      &
		  	\psfrag{valf}[b][t]{}
			\psfrag{xlabel}[t][b]{CPU Time (sec.)}
			\psfrag{Grad}{\scriptsize Gradient}
			\psfrag{GradNest}{\scriptsize Nesterov}
			\psfrag{GradNestBack}{\scriptsize Nest. + back.}
			\psfrag{RMPE5}{\scriptsize RNA 5}
			\psfrag{RMPE5LS}{\scriptsize RNA 5 + LS}
		    \includegraphics[width=0.45\textwidth]{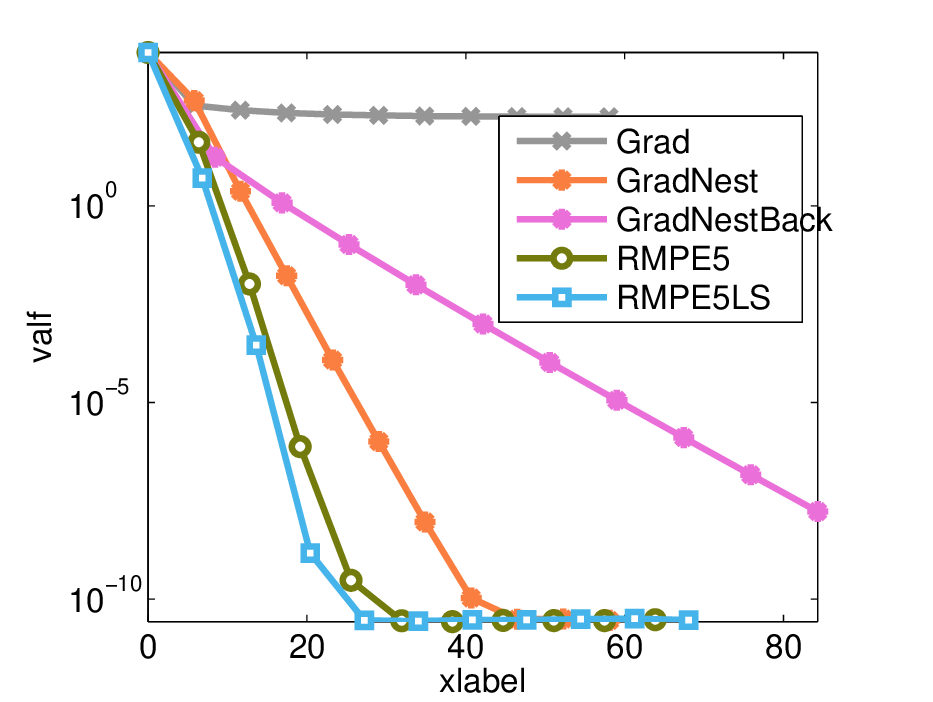}
		\end{tabular}
	\end{center}
	   \caption{ Logistic regression on sido0 dataset, with $\tau = 10^{2}$ (condition number = $1.5\cdot 10^{5}$) }
	   \label{fig:sido0}
\end{figure}

\begin{figure}[h]
	\begin{center}
		\begin{tabular}{cc}
		  	\psfrag{valf}[b][t]{$f(x_k)-f(x^*)$}
			\psfrag{xlabel}[t][b]{Gradient oracle calls}
			\psfrag{Grad}{\scriptsize Gradient}
			\psfrag{GradNest}{\scriptsize Nesterov}
			\psfrag{GradNestBack}{\scriptsize Nest. + back.}
			\psfrag{RMPE5}{\scriptsize RNA 5}
			\psfrag{RMPE5LS}{\scriptsize RNA 5 + LS}
			\includegraphics[width=0.48\textwidth]{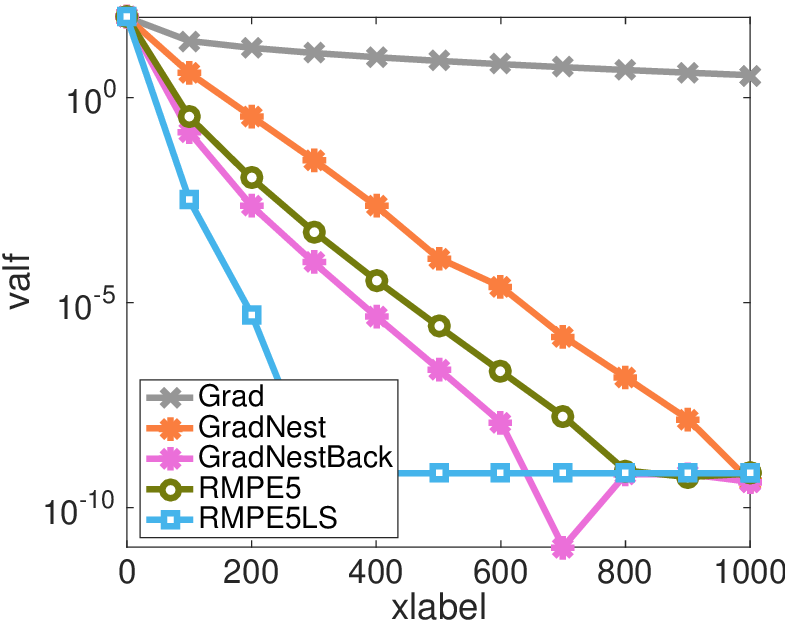}
		      &
		  	\psfrag{valf}[b][t]{}
			\psfrag{xlabel}[t][b]{CPU Time (sec.)}
			\psfrag{Grad}{\scriptsize Gradient}
			\psfrag{GradNest}{\scriptsize Nesterov}
			\psfrag{GradNestBack}{\scriptsize Nest. + back.}
			\psfrag{RMPE5}{\scriptsize RNA 5}
			\psfrag{RMPE5LS}{\scriptsize RNA 5 + LS}
		    \includegraphics[width=0.45\textwidth]{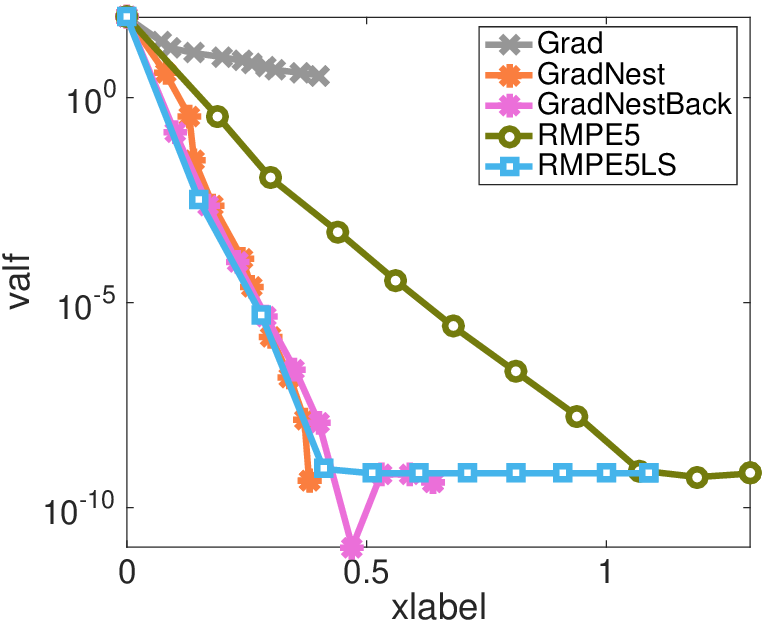}\\
		    \psfrag{valf}[b][t]{$f(x_k)-f(x^*)$}
   			\psfrag{xlabel}[t][b]{Gradient oracle calls}
   			\psfrag{Grad}{\scriptsize Gradient}
   			\psfrag{GradNest}{\scriptsize Nesterov}
   			\psfrag{GradNestBack}{\scriptsize Nest. + back.}
   			\psfrag{RMPE5}{\scriptsize RNA 5}
   			\psfrag{RMPE5LS}{\scriptsize RNA 5 + LS}
   			\includegraphics[width=0.45\textwidth]{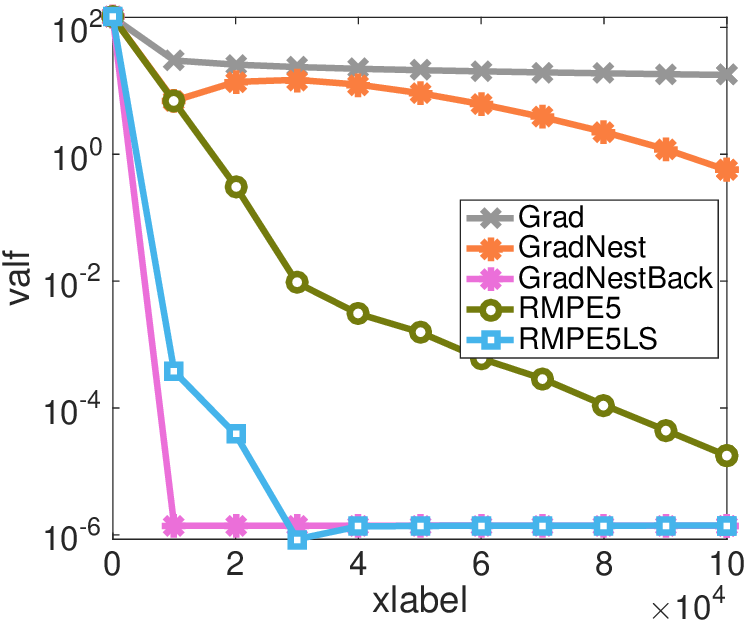}
   			&
   			\psfrag{valf}[b][t]{}
   			\psfrag{xlabel}[t][b]{CPU Time (sec.)}
   			\psfrag{Grad}{\scriptsize Gradient}
   			\psfrag{GradNest}{\scriptsize Nesterov}
   			\psfrag{GradNestBack}{\scriptsize Nest. + back.}
   			\psfrag{RMPE5}{\scriptsize RNA 5}
   			\psfrag{RMPE5LS}{\scriptsize RNA 5 + LS}
   			\includegraphics[width=0.45\textwidth]{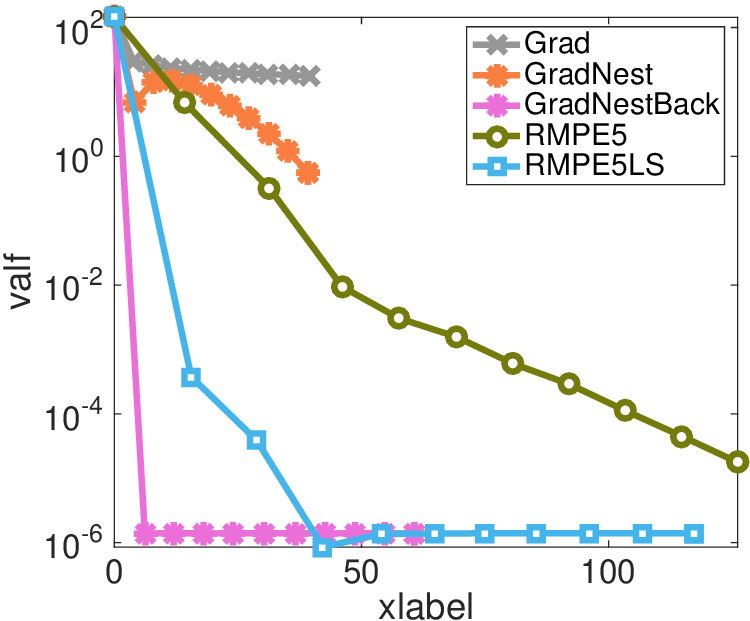}
		\end{tabular}
	\end{center}
	   \caption{ Logistic regression on sonar dataset. From top to bottom, we used $\tau = 10^{-1}$ (condition number = $7\cdot 10^{3}$) and $\tau = 10^{-6}$ (condition number = $7\cdot 10^{8}$).}
	   \label{fig:sonar}
\end{figure}

\begin{figure}[h]
	\begin{center}
		\begin{tabular}{cc}
		  	\psfrag{valf}[b][t]{$f(x_k)-f(x^*)$}
			\psfrag{xlabel}[t][b]{Gradient oracle calls}
			\psfrag{Grad}{\scriptsize Gradient}
			\psfrag{GradNest}{\scriptsize Nesterov}
			\psfrag{GradNestBack}{\scriptsize Nest. + back.}
			\psfrag{RMPE5}{\scriptsize RNA 5}
			\psfrag{RMPE5LS}{\scriptsize RNA 5 + LS}
			\includegraphics[width=0.45\textwidth]{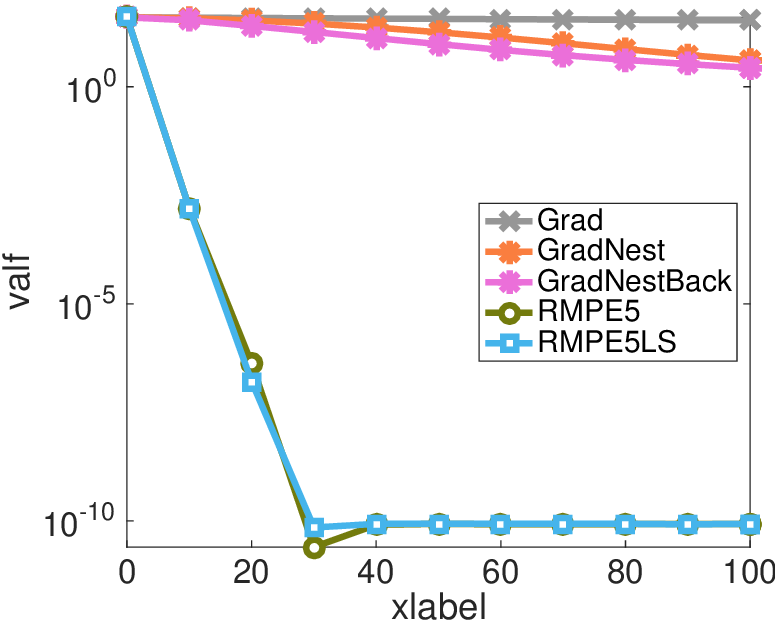}
		      &
		  	\psfrag{valf}[b][t]{}
			\psfrag{xlabel}[t][b]{CPU Time (sec.)}
			\psfrag{Grad}{\scriptsize Gradient}
			\psfrag{GradNest}{\scriptsize Nesterov}
			\psfrag{GradNestBack}{\scriptsize Nest. + back.}
			\psfrag{RMPE5}{\scriptsize RNA 5}
			\psfrag{RMPE5LS}{\scriptsize RNA 5 + LS}
		    \includegraphics[width=0.45\textwidth]{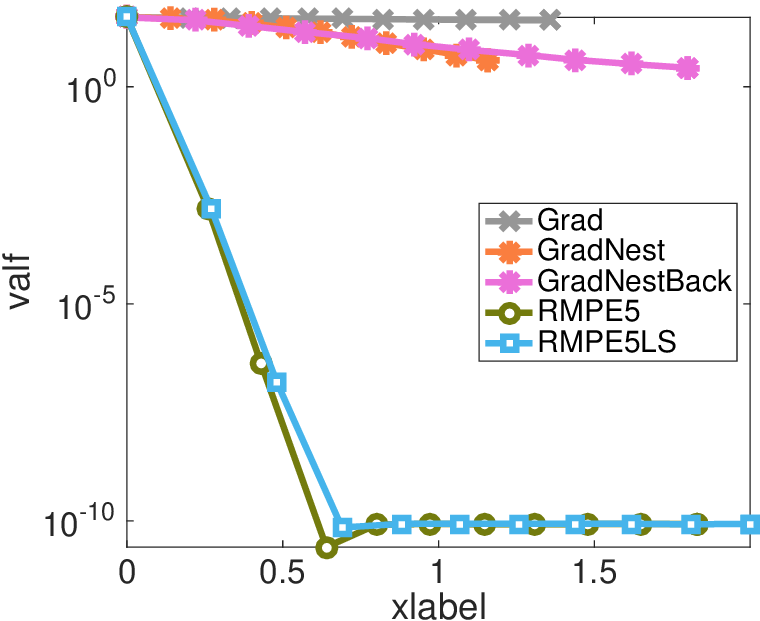} \\
		    \psfrag{valf}[b][t]{$f(x_k)-f(x^*)$}
   			\psfrag{xlabel}[t][b]{Gradient oracle calls}
   			\psfrag{Grad}{\scriptsize Gradient}
   			\psfrag{GradNest}{\scriptsize Nesterov}
   			\psfrag{GradNestBack}{\scriptsize Nest. + back.}
   			\psfrag{RMPE5}{\scriptsize RNA 5}
   			\psfrag{RMPE5LS}{\scriptsize RNA 5 + LS}
   			\includegraphics[width=0.45\textwidth]{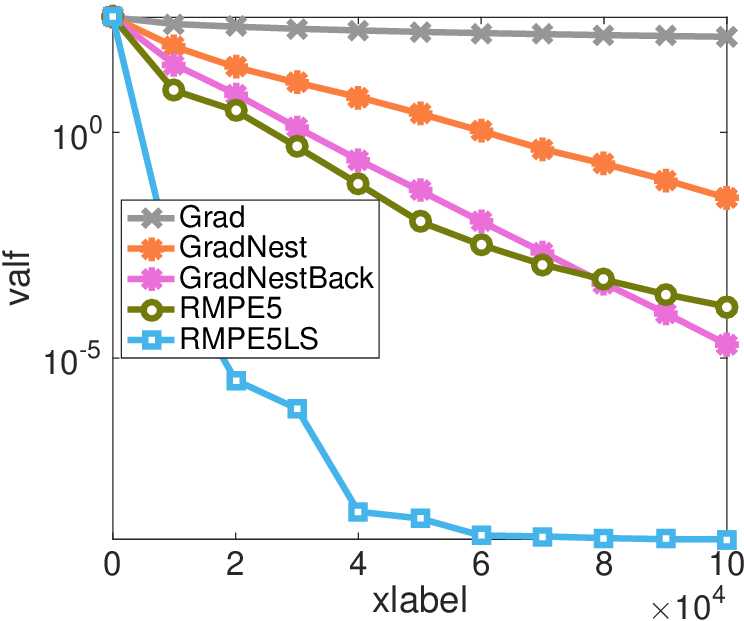}
  			&
   		  	\psfrag{valf}[b][t]{}
   			\psfrag{xlabel}[t][b]{CPU Time (sec.)}
   			\psfrag{Grad}{\scriptsize Gradient}
   			\psfrag{GradNest}{\scriptsize Nesterov}
   			\psfrag{GradNestBack}{\scriptsize Nest. + back.}
   			\psfrag{RMPE5}{\scriptsize RNA 5}
   			\psfrag{RMPE5LS}{\scriptsize RNA 5 + LS}
			\includegraphics[width=0.45\textwidth]{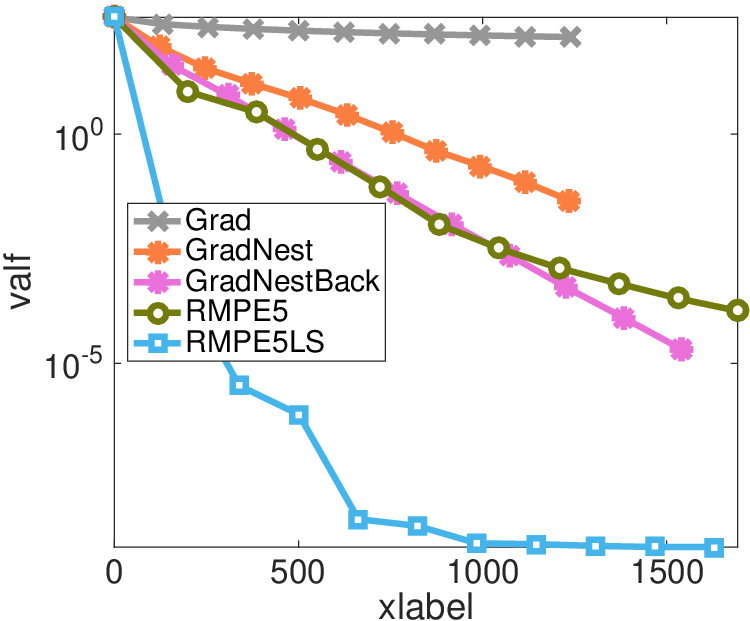}\\
			\psfrag{valf}[b][t]{$f(x_k)-f(x^*)$}
			\psfrag{xlabel}[t][b]{Gradient oracle calls}
			\psfrag{Grad}{\scriptsize Gradient}
			\psfrag{GradNest}{\scriptsize Nesterov}
			\psfrag{GradNestBack}{\scriptsize Nest. + back.}
			\psfrag{RMPE5}{\scriptsize RNA 5}
			\psfrag{RMPE5LS}{\scriptsize RNA 5 + LS}
			\includegraphics[width=0.45\textwidth]{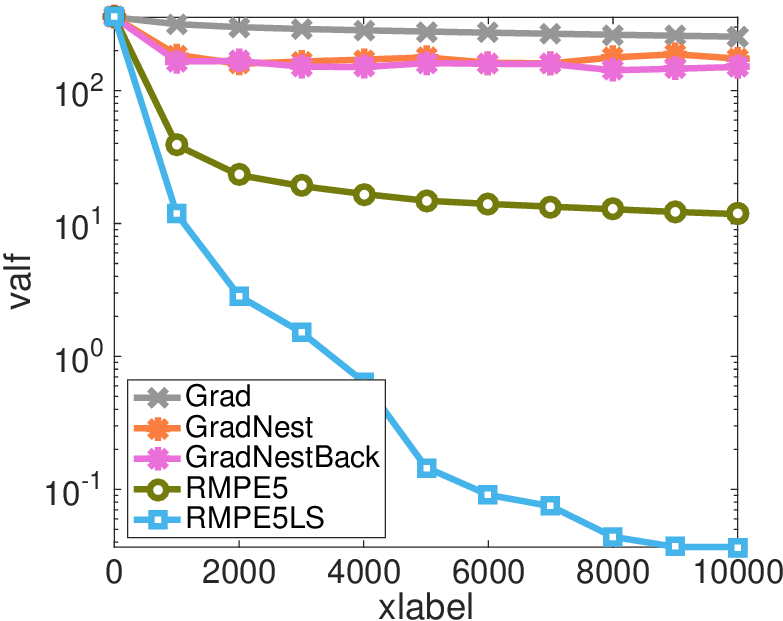}
			&
			\psfrag{valf}[b][t]{}
			\psfrag{xlabel}[t][b]{CPU Time (sec.)}
			\psfrag{Grad}{\scriptsize Gradient}
			\psfrag{GradNest}{\scriptsize Nesterov}
			\psfrag{GradNestBack}{\scriptsize Nest. + back.}
			\psfrag{RMPE5}{\scriptsize RNA 5}
			\psfrag{RMPE5LS}{\scriptsize RNA 5 + LS}
			\includegraphics[width=0.45\textwidth]{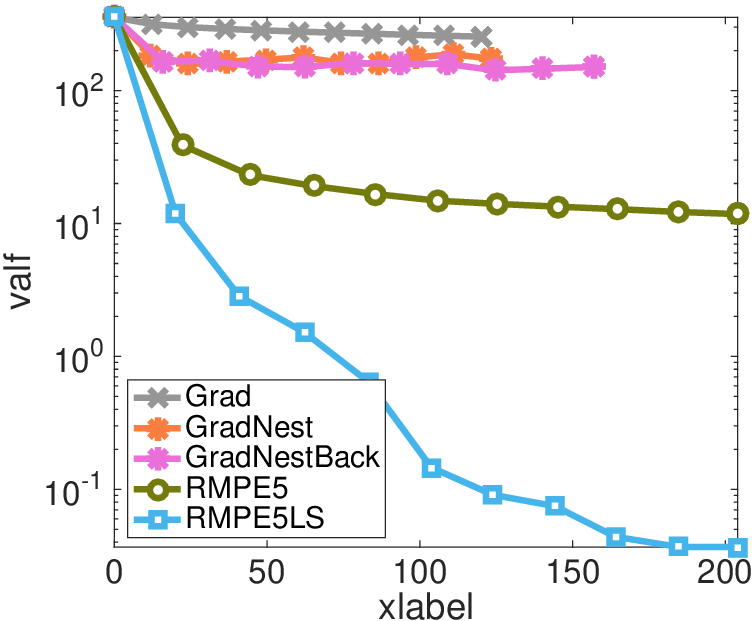}
		\end{tabular}
	\end{center}
   \caption{Logistic regression on Madelon dataset. From top to bottom, we used , $\tau = 10^7$ (condition number = $6\cdot 10^3$), $\tau = 10^2$ (condition number = $1.2\cdot10^9$ and $\tau = 10^{-3}$ (condition number = $6\cdot 10^{13}$).}
   \label{fig:madelon}
\end{figure}

Figure \ref{fig:unstableAMPE} shows the importance of the regularization in the acceleration algorithm. Indeed, if we use Algorithm \ref{algo:acc_fixedpoint} then the norm of the inverse of $\tilde R^T \tilde R$ may be huge, so the computation of the coefficients $\tilde c_\lambda^*$ is unstable. This leads to an unreliable acceleration method, which may improve sometimes the accuracy, but often making the process divergent. In Figures \ref{fig:sido0}, \ref{fig:sonar} and \ref{fig:madelon}, we see that our algorithm has a similar behavior to the conjugate gradient: unlike the Nesterov's method, where we need to provide  parameters $\mu$ and $L$, the acceleration algorithm adapts himself in function of the spectrum of $G$ (so it can exploit the good local strong convexity parameter), without any prior specification. We can, for example, observe this behavior when the global strong convexity parameter is bad but not the local one.

\section{Conclusion and Perspectives}
In this paper, we developed a method which is able to accelerate, under some regularity conditions, the convergence of a sequence $\{\tilde x_i\}$ without any information on the algorithm which generated this sequence. The regularization parameter used in the acceleration method is found by a simple and inexpensive grid-search. The algorithm itself is simple as it only requires solving a small linear system. Also, we showed (using gradient method on logistic regression) that the strategy which consists in restarting the algorithm after an extrapolation method can lead to significantly improved convergence rates.
Future work will consist in improving the performance of the algorithm by exploiting the structure of the perturbations matrix in some cases and extending the algorithm to the stochastic case and to the non-symmetric case.

\clearpage
{
\bibliographystyle{agsm}
\bibliography{reg_non_acc}
}

\section*{Acknowledgements}
AA is at the d\'epartement d'informatique de l'ENS, \'Ecole normale sup\'erieure, UMR CNRS 8548, PSL Research University, 75005 Paris, France, and INRIA Sierra project-team. The authors would like to acknowledge support from a starting grant from the European Research Council (ERC project SIPA), from the ITN MacSeNet (project number 642685), as well as support from the chaire {\em \'Economie des nouvelles donn\'ees} with the {\em data science} joint research initiative with the {\em fonds AXA pour la recherche}, and from a Google focused award.

\clearpage

\appendix

\section{Missing propositions and proofs}
\label{ss:missongproofs}
\subsection{Missing propositions}

\begin{proposition} \label{prop:opt_val_sqrt_fun}
Consider the function
\[
	f(x) = \sqrt{a - \lambda x^2} + bx 
\]
defined for $x \in [0,\sqrt{a/\lambda}]$. The its maximal value is attained at
\[
	x_{ \text{opt} } = \frac{b \sqrt{a}}{\sqrt{\lambda^2\kappa^2 + \lambda b^2}}
\]
and its maximal value is thus, if $x_{ \text{opt} } \in [0,\sqrt{a/\lambda}]$,
\BEA
	f_{\max} =  \sqrt{a} \sqrt{\kappa^2 + \frac{b^2}{\lambda}} \label{eq:opt_val_sqrt_fun}.
\EEA
\end{proposition}
\begin{proof}
The (positive) root of the derivative of $f$ follows
\[
	b\sqrt{a-\lambda x^2} - \kappa\lambda x = 0 \qquad \Leftrightarrow \qquad x = \frac{b \sqrt{a}}{\sqrt{\lambda^2\kappa^2 + \lambda b^2}} .
\]
If we inject the solution in our function, we obtain its maximal value,
\BEAS
	\kappa \sqrt{a - \lambda \left(\frac{b \sqrt{a}}{\sqrt{\lambda^2\kappa^2 + \lambda b^2}}\right)^2} + b \frac{b \sqrt{a}}{\sqrt{\lambda^2\kappa^2 + \lambda b^2}} & = & \kappa \sqrt{a - \lambda \frac{b^2 a}{\lambda^2\kappa^2 + \lambda b^2}} + b \frac{b \sqrt{a}}{\sqrt{\lambda^2\kappa^2 + \lambda b^2}}, \\
	& = & \kappa \sqrt{a - \lambda \frac{b^2 a}{\lambda^2\kappa^2 + \lambda b^2}} + b \frac{b \sqrt{a}}{\sqrt{\lambda^2\kappa^2 + \lambda b^2}}, \\
	& = & \kappa \sqrt{\frac{a \lambda^2 \kappa^2}{\lambda^2\kappa^2 + \lambda b^2}} + b \frac{b \sqrt{a}}{\sqrt{\lambda^2\kappa^2 + \lambda b^2}}, \\
	& = &  \sqrt{a} \frac{ \kappa^2 \lambda + b^2 }{\sqrt{\lambda^2\kappa^2 + \lambda b^2}}, \\
	& = &  \frac{\sqrt{a}}{\lambda} \sqrt{\lambda^2\kappa^2 + \lambda b^2}.
\EEAS
The simplification with $\lambda$ in the last equality concludes the proof.
\end{proof}

\subsection{Proof of proposition \ref{prop:explicit_bounds_gradient}} \label{app:proof_explicit_bounds_gradient}
First, we show that the choice $\sigma = 1-\frac{\mu}{L}$ satisfies $\|G\| = \|g'(x^*)\| \leq \sigma$. Our fixed-point function $g$ reads
\[
	g(x) = x-\frac{1}{L}f'(x) .
\]
Since $g'(x) = I-\frac{1}{L}f''(x)$, we have $g'(x^*) = I-\frac{1}{L}f''(x^*)$. Because $f$ is $\mu$-strongly convex, $f''(x)\succeq \mu I$, in particular at $x=x^*$. In conclusion,
\[
	\|g'(x^*)\| = \|I-\frac{1}{L}f''(x^*)\| \leq 1-\frac{\mu}{L}.
\]
Now, consider the matrix $\tilde R$. Since the $i-th$ column $\tilde R_i$ is equal to $\tilde x_{i}-\tilde x_{i-1}$,
\BEAS
	\|\tilde R_i\| & = & \|\tilde x_{i}-\tilde x_{i-1}\|, \\
	& = & \frac{1}{L}\|f'(\tilde x_i)\|, \\
	& \leq & \|\tilde x_i-x^*\| .
\EEAS
In the last inequality, we used the fact that $f$ is $L$-Lipschitz, so $\|f(x)-f(x^*)\| \leq L \|x-x^*\|$. It is also possible to prove \citep{Nest03a} that gradient method converges at rate
\[
	\|\tilde x_{i+1}-x^*\| \leq \sigma \|x_i-x^*\|.
\]
It remains to link this quantity to $\|\tilde R\|$,
\BEAS
	\|\tilde R\| & \leq & \sum_{i=0}^k \|R_i\|, \\
	& \leq & \sum_{i=0}^k \sigma^i\|x_0-x^*\|, \\
	& = & \frac{1-\sigma^{k+1}}{1-\sigma}\|x_0-x^*\|. \\
\EEAS
We continue with $\|\mathcal{E}\|$. We express $\|\mathcal{E}_i\| = \|\tilde x_{i+1}-x_{i+1}\|_2$ in function of $\|\tilde x_0-x_0\|_2$ using a recursion with $\|\tilde x_i-x_i\|_2$,
\BEAS\tilde x_{i+1} - x_{i+1} & = & \tilde x_i - \frac{1}{L}\nabla f(\tilde x_i) - x_i + \frac{1}{L} \nabla^2 f(x^*)(x_i-x^*) , \\
& = & \tilde x_{i} - x_{i} - \frac{1}{L}(\nabla f(\tilde x_i) - \nabla^2 f(x^*)(x_i-x^*)) , \\
& = & \left(I-\frac{\nabla^2 f(x^*)}{L}\right)(\tilde x_{i} - x_{i}) - \frac{1}{L}(\nabla f(\tilde x_i) - \nabla^2 f(x^*)(\tilde x_i-x^*)) .
\EEAS
Since our function has a Lipschitz-continuous Hessian, it is possible to show that (\cite{Nest03a}, \mbox{Lemma 1.2.4})
\BEQ
\left\|\nabla f(y) - \nabla f(x) - \nabla^2 f(x)(y-x)\right\|_2 \leq \frac{M}{2}\|y-x\|^2.
\label{eq:boundhessianlipch}
\EEQ
We can thus bound the norm of the error at the $i^{\text{th}}$ iteration,
\BEAS
\|x_{i+1}-\tilde x_{i+1}\|_2 & \leq & \left\| I-\frac{\nabla^2 f(x^*)}{L})\right\|_2 \| x_{i}-\tilde x_{i} \|_2 + \frac{1}{L}\left\|\nabla f(\tilde x_i) - \nabla^2 f(x^*)(\tilde x_i-x^*)\right\|_2,\\
& = & \| g''(x^*)\|_2\| x_{i}-\tilde x_{i} \|_2 + \frac{1}{L}\left\|\nabla f(\tilde x_i) - \nabla f(x^*) - \nabla^2 f(x^*)(\tilde x_i-x^*)\right\|_2.\\
\EEAS
By equation \eqref{eq:boundhessianlipch}, and because $\|g''(x^*)\| \leq \sigma$, we have
\BEAS
\|x_{i+1}-\tilde x_{i+1}\|_2 & \leq & \sigma \|  x_{i}-\tilde x_{i} \|_2 + \frac{M}{2L}\left\|\tilde x_i - x^*\right\|_2^2 , \\
& \leq & \sigma \| x_{i}-\tilde x_{i}\|_2 + \frac{M}{2L}\sigma^{2i} \|x_0-x^*\|_2^2 , \\
& \leq & \| x_{i}-\tilde x_{i}\|_2 + \frac{M}{2L} \|x_0-x^*\|_2^2.
\EEAS
The simplification in the last line greatly simplifies future computations. We thus have the bound
\[
	\|x_{i+1}-\tilde x_{i+1}\|_2 \leq (i+1) \frac{M}{2L}\|x_0-x^*\|^2.
\]
Finally,
\BEAS
	\|\mathcal{E}\| & \leq & \sum_{i=0}^k \|x_{i+1}-\tilde x_{i+1}\|_2 , \\
	& \leq & \sum_{i=0}^k (i+1) \frac{M}{2L}\|x_0-x^*\|^2 , \\
	& \leq & (k+2)^2 \frac{M}{4L}\|x_0-x^*\|^2.
\EEAS
Despite the simplification made earlier, the results of this bounds are close to the one obtained without simplification.

\clearpage
\end{document}